%% file: locallygrid_MAIN_FILE_-_edited.tex
\documentclass[reqno,11pt,a4paper]{amsart}

\usepackage{amsmath,amssymb,amsthm,enumerate,graphicx,mathrsfs,pstricks,pst-node,multido}

\newcommand{\V}{\textnormal{V}}
\newcommand{\E}{\textnormal{E}}
\newcommand{\val}{\textnormal{val}}
\newcommand{\diam}{\textnormal{diam}}
\newcommand{\F}{\mathbb{F}}
\newcommand{\Aut}{\textnormal{Aut}}
\newcommand{\Sp}{\textnormal{Sp}}

\newtheorem{theorem}{Theorem}[section]
\newtheorem{lemma}[theorem]{Lemma}
\newtheorem{corollary}[theorem]{Corollary}
\newtheorem{proposition}[theorem]{Proposition}

\theoremstyle{definition}

\newtheorem{remark}[theorem]{Remark}

\newtheorem{construction}[theorem]{Construction}
\newtheorem*{problem}{Problem}

\numberwithin{equation}{section}

\def\mucomponent{
	\psline[linewidth=1.4pt](-1,1)(-0.5,1)(-0.5,0.5)(0,0.5)(0,0)
	\psline[linestyle=dotted,linewidth=1.4pt](0,0)(0.5,-0.5)
	\psline[linewidth=1.4pt](0.5,-0.5)(1,-0.5)(1,-1)
	\pnode(1,-1){a} \pnode(-1,-1){b} \pnode(-1,1){c}
	\ncarc[arcangle=15,linewidth=1.4pt]{a}{b} \ncarc[arcangle=15,linewidth=1.4pt]{b}{c}
}
\def\localgrid{
\pspolygon[linecolor=lightgray](-1.5,-1.5)(-1.5,1.5)(1.5,1.5)(1.5,-1.5)
\pspolygon[linecolor=lightgray,fillstyle=hlines,hatchcolor=lightgray,hatchangle=0](0,1.5)(1.5,1.5)(1.5,0)(0,0)(0,1.5)
\pspolygon[linecolor=lightgray,fillstyle=vlines,hatchcolor=lightgray,hatchangle=0](0,-1.5)(-1.5,-1.5)(-1.5,0)(0,0)(0,-1.5)
}
\def\distancediagram{
	\rput(-7,0){\ovalnode{k0}{$1$}} \rput(-6.4,0.3){\small $n^2$}
	\rput(-5,0){\ovalnode{k1}{$n^2$}} \rput(-5.6,0.25){\small $1$} \rput(-3.9,0.25){\small $(n-1)^2$}
	\rput(-1.4,0){\ovalnode{k2}{$k_2(x)$}}
	\rput(4,0){\ovalnode{k3}{$k_3(x)$}}
	\pnode(7.2,0){A}
	\ncline{k0}{k1} \ncline{k1}{k2} \ncline{k2}{k3} \ncline{k3}{A}
	\rput(7.7,0){$\ldots$}
	}
\def\distancediagramsmall{
	\rput(0,0){\ovalnode{ki}{$k_i(x)$}} \pnode(-2.5,0){A} \pnode(3.1,0){B}
	\ncline{A}{ki} \ncline{ki}{B} \rput(-3,0){$\ldots$} \rput(3.6,0){$\ldots$}
	\rput(0,-0.6){\footnotesize ($i \geq 4$)}
	}
\def\xdiagram{
	\pscircle*[linecolor=lightgray](-1.5,0.75){0.4cm} \psellipse*[linecolor=lightgray](1.5,0)(0.55,0.65)
	\pnode(-1.5,0.75){e0} \pnode(1.5,0){C} \ncline[linecolor=lightgray,linewidth=0.25cm]{C}{e0}
	\rput(-1.5,0){\psellipse(0,0)(1,1.5)} \rput(1.5,0){\psellipse(0,0)(1,1.5)} 
	\rput(-1.5,1.8){\small $\Gamma(x)$} \rput(1.5,1.8){\small $\Gamma_2(x)$}
	}
\def\gridgraph3{
	\multido{\n=-1+1}{3}{\psline[linecolor=lightgray](\n,-1)(\n,1)}
	\multido{\n=-1+1}{3}{\psline[linecolor=lightgray](-1,\n)(1,\n)}
	\rput(0,-0.732){\multido{\n=-4+1}{3}{\psarc[linecolor=lightgray](0,\n){3.864}{75}{105}}}
	\rput{90}(0.732,0){\multido{\n=-4+1}{3}{\psarc[linecolor=lightgray](0,\n){3.864}{75}{105}}}
	}

\title{On locally ${n} \times {n}$ grid graphs \today}
\date{\today}
\author{Carmen Amarra} \address{Institute of Mathematics, University of the Philippines Diliman}
\author{Wei Jin} \address{School of Statistics, Jiangxi University of Finance and Economics}
\author{Cheryl E. Praeger} \address{Centre for the Mathematics of Symmetry and Computation, The University of Western Australia}

\begin{document}

\begin{abstract}
We investigate locally $n \times n$ grid graphs, that is, graphs in which the neighbourhood of any vertex is the Cartesian product of two complete graphs on $n$ vertices. We consider the subclass of these graphs for which each pair of vertices at distance two is joined by sufficiently many paths of length $2$. The number of such paths is known to be at most $2n$ by previous work of Blokhuis and Brouwer. We show that 
if each pair is joined by at least $2(n-1)$ such paths then the diameter is at most $3$ and we give a tight upper bound on the order of the graphs. We show that graphs meeting this upper bound are distance-regular antipodal covers of complete graphs. We exhibit an infinite family of such graphs which are locally $n \times n$ grid for odd prime powers $n$, and apply these results to locally $5 \times 5$ grid graphs to obtain a classification for the case where either all $\mu$-graphs have order at least $8$ or all $\mu$-graphs have order $c$ for some constant $c$.
\end{abstract}

\maketitle

\input{locallygrid-sec1-intro}
\input{locallygrid-sec2-prelims-sec3-examples}
\input{locallygrid-sec4-basics}
\input{locallygrid-sec5-cliques}
\input{locallygrid-sec6-lowerbound}
\input{locallygrid-sec7-appendix}

\end{document}

%% file: locallygrid-sec1-intro.tex
\section{Introduction}

Throughout this paper all graphs are finite, simple, and undirected.

Let $m$ and $n$ be integers. An \emph{$m \times n$ grid} (also known as the \emph{$m \times n$ lattice graph}) is the Cartesian product $K_m \square K_n$ of two complete graphs, one with order $m$ and the other with order $n$. It has as vertices all ordered pairs $(i,j)$, $i \in \{1, \ldots, m\}$ and $j \in \{1, \ldots, n\}$, and as edges all $2$-sets of ordered pairs that agree in exactly one coordinate. If $m, n \geq 2$ then an $m \times n$ grid has diameter $2$. The $n \times n$ grid, sometimes called the lattice graph $L_2(n)$ of order $n$, is isomorphic to the Hamming graph $H(2,n)$, and has automorphism group $S_n \wr S_2$ which acts transitively of rank $3$ on its vertex set.

For any class $\mathcal{G}$ of graphs, a graph is said to be \emph{locally $\mathcal{G}$} if the induced subgraph on the neighbourhood of any vertex is isomorphic to a graph in $\mathcal{G}$. In particular, a graph is said to be \emph{locally grid} if $\mathcal{G}$ is the class of all grid graphs. Locally grid graphs were first studied in 1977 by Buekenhout and Hubaut in \cite{localpolar}, where they arise as adjacency graphs of certain locally polar spaces. In particular, they exhibit two infinite families of graphs which provide examples of locally $n \times n$ grid graphs for all $n \geq 4$ \cite[Section 2.3]{localpolar}. (One of these families consists of the Johnson graphs, and the other consists of quotients of the Johnson graphs by an antipodal partition. For $n \geq m$ the Johnson graph $J(n,m)$ has vertices the set of all $m$-subsets of an $n$-set, and two $m$-subsets are adjacent if their intersection has size $m-1$. If $n=2m$, the quotient of $J(2m,m)$ by an antipodal partition is the \emph{half Johnson graph} $\frac{1}{2}J(2m,m)$.) Families of locally $m \times n$ grid graphs which have been completely classified include the subcases where $m = 2$ \cite{4by4}, $m = 3$ \cite{3byq} (see also Remark \ref{remark:3x3}), and $m = n = 4$ \cite{4by4}. The first are all triangular graphs, and the second are line graphs of certain connected partial linear spaces. The third classification, for $m = n = 4$, yields exactly four graphs, namely the Johnson graph $J(8,4)$ and its quotient $\frac{1}{2}J(8,4)$, and two graphs on $40$ vertices.


A \emph{$\mu$-graph} of a non-complete graph is an induced subgraph on the set of common neighbours of two vertices at distance two. Blokhuis and Brouwer showed in \cite{4by4} that any $\mu$-graph of a locally grid graph is a union of cycles of even length, and that if each $\mu$-graph is a union of $4$-cycles then the graph is either a Johnson graph or a quotient of a Johnson graph. Furthermore, if all $\mu$-graphs have the maximum possible order (that is, $2m$ if the graph is locally $m \times n$ grid with $m \leq n$) then $\Gamma$ is strongly regular and the parameters are known. In \cite{bilinear}, Gavrilyuk and Koolen considered locally $m \times n$ grid graphs with $n \geq m \geq 3$ whose $\mu$-graphs are all $6$-cycles, and with the additional condition that for each pair of vertices $x$ and $y$ at distance two, there are $(m-3)(n-3)$ vertices adjacent to $y$ and at distance three from $x$. They characterised such graphs as certain quotients of the graph of bilinear $(d \times e)$-forms over the field $\F_2$ where $m = 2^d-1$ and $n = 2^e-1$.

In this paper we undertake a general study of locally $n \times n$ grid graphs extending some of the results in \cite{4by4}. Our first result is a general characterisation for the case where all $\mu$-graphs are nearly as large as possible. We denote the vertex set of the graph $\Gamma$ by $\V(\Gamma)$ and its diameter by $\diam(\Gamma)$.

\begin{theorem} \label{maintheorem:order-diameter}
Assume that $\Gamma$ is connected and locally $n \times n$ grid for some $n \geq 2$. Then any $\mu$-graph has even order at least $4$ and at most $2n$.
	\begin{enumerate}[(1)]
	\item \label{>2(n-1)} If all $\mu$-graphs of $\Gamma$ have order at least $2(n-1)$, then
		\[ |\V(\Gamma)| \leq \left\lfloor \frac{(n^2 + 1)(n+1)}{2} \right\rfloor \quad \text{and} \quad \diam(\Gamma) \leq 3. \]
	\item \label{=2(n-1)} In part (\ref{>2(n-1)}), $|\V(\Gamma)| = \lfloor (n^2 + 1)(n+1)/2 \rfloor$ if and only if all $\mu$-graphs have order equal to $2(n-1)$, and in this case $n$ is odd, $\diam(\Gamma) = 3$, and $\Gamma$ is a distance-regular antipodal $((n+1)/2)$-cover of $K_{n^2+1}$ with intersection array $\left(n^2,\, (n-1)^2,\, 1;\, 1,\, 2(n-1),\, n^2\right)$.
	\end{enumerate}
\end{theorem}


There are examples of graphs satisfying the conditions of Theorem \ref{maintheorem:order-diameter} (\ref{=2(n-1)}). An infinite family of such graphs arises from a construction of Godsil and Hensel \cite{DRcovers}, which in turn is a special case of the construction given in \cite[Proposition 12.5.3]{BCN}. We describe this in Construction \ref{example:drg}.

\begin{theorem} \label{maintheorem:family}
For each odd prime power $n$ the graph $\Gamma^{(n)}$ in Construction \ref{example:drg} is locally $n \times n$ grid and satisfies the conditions in Theorem \ref{maintheorem:order-diameter} (\ref{=2(n-1)}). Furthermore each $\mu$-graph in $\Gamma^{(n)}$ is either connected or a union of cycles of equal length.
\end{theorem}

Theorem \ref{maintheorem:family} will follow from the technical Proposition \ref{proposition:family} which gives, in addition, local structural information and describes the $\mu$-graphs for the graphs in Construction \ref{example:drg}. In particular we show that the number of cycles in a $\mu$-graph is unbounded (see Proposition \ref{proposition:family} (\ref{cycles-lowerbound})).

In addition to the above, we also obtain technical results about maximal cliques in general locally $n \times n$ grid graphs. We apply these together with Theorem \ref{maintheorem:order-diameter} to the case where $n = 5$, and obtain the following.

\begin{theorem} \label{maintheorem:5by5}
Assume that $\Gamma$ is connected and locally $5 \times 5$ grid.
	\begin{enumerate}[(1)]
	\item \label{mu>8} If all $\mu$-graphs in $\Gamma$ have order at least $8$, then $\Gamma$ has diameter $3$. Moreover, either
            \begin{enumerate}[(i)]
            \item \label{thm:order72} $\Gamma$ has $72$ vertices, some $\mu$-graphs of $\Gamma$ have order $8$ and some have order $10$, $\Gamma$ is an antipodal double cover of $K_6 \times K_6$, and has distance diagram as in Figure \ref{figure:diagram-5x5}, or
            \item \label{thm:order78} all $\mu$-graphs of $\Gamma$ have order equal to $8$, and $\Gamma$ is a distance-regular antipodal triple cover of $K_{26}$ with intersection array $(25, 16, 1; 1, 8, 25)$.
            \end{enumerate}
	\item \label{mu-constant} If all $\mu$-graphs in $\Gamma$ have constant order $|\mu|$, then either $|\mu| = 8$ and $\Gamma$ is as in part \eqref{mu>8} {\rm (ii)}, or $|\mu| = 4$ and $\Gamma$ is the Johnson graph $J(10,5)$ or the half Johnson graph $\frac{1}{2}J(10,5)$.
	\end{enumerate}
\end{theorem}

\begin{figure}[h]
\begin{center}
\begin{pspicture}(-4.5,-2.5)(4.5,0.7)
\rput(-4.5,0){\rput(0,0){\ovalnode{k0}{$1$}} \rput(0.6,0.25){\small $25$}}
\rput(-1.5,0){\rput(0,0){\ovalnode{k1}{$25$}} \rput(-0.6,0.25){\small $1$} \rput(0,0.55){\small $8$} \rput(0.6,0.25){\small $8$} \rput(0,-0.55){\small $8$}}
\rput(0,-2){\rput(0,0){\ovalnode{k210}{$20$}} \rput(0,-0.55){\small $5$} \rput(-0.65,0.35){\small $10$} \rput(0.65,0.35){\small $10$}}
\rput(1.5,0){\rput(0,0){\ovalnode{k28}{$25$}} \rput(0.6,0.25){\small $1$} \rput(0,0.55){\small $8$} \rput(-0.6,0.25){\small $8$} \rput(0,-0.55){\small $8$}}
\rput(4.5,0){\rput(0,0){\ovalnode{k3}{$1$}} \rput(-0.6,0.25){\small $25$}}
\ncline{k0}{k1} \ncline{k1}{k28} \ncline{k28}{k3} \ncline{k1}{k210} \ncline{k28}{k210}
\end{pspicture}
\end{center}
\caption{Distance diagram for $\Gamma$ in Theorem \ref{maintheorem:5by5} (\ref{order72})}
\label{figure:diagram-5x5}
\end{figure}

There is at least one graph satisfying the conditions of Theorem \ref{maintheorem:5by5} (\ref{mu>8}) with all $\mu$-graphs of order $8$, namely, the graph in Construction \ref{example:drg} with $n = 5$. Likewise, a graph satisfying the conditions of Theorem \ref{maintheorem:5by5} \eqref{mu>8} with $|\V(\Gamma)| = 72$ exists; this is described in \cite[Example 9.15(iii)]{egq2} as arising from a $2$-fold cover of an extended generalised quadrangle (see Subsection \ref{ss:egq}) EGQ(4,1), and in \cite[Proposition 3.1]{egq} as arising from a hyperoval in the rank $3$ polar space $Q^+(5,4)$. Dima Pasechnik has given additional information about the geometric structure of this example in \cite{ho}. In this paper he also constructs a locally $5 \times 5$ grid graph of order $96$ corresponding to a hyperoval in $Q^+(6,4)$. It has $\mu$-graphs of orders $4$ and $8$. It would be interesting to know all the locally $5 \times 5$ grids. (By Lemma 6.4, the orders of any new examples are at most $270$.)

\begin{problem}
Complete the classification of connected locally $5 \times 5$ grid graphs.
\end{problem}

The rest of the paper is organised as follows: In Section \ref{sec:prelims} we list elementary properties of locally $m \times n$ grid graphs. In Section \ref{sec:family} we introduce the infinite family of graphs mentioned above, and prove Theorem \ref{maintheorem:family}. We then restrict ourselves to the case where $m = n$, and in Section \ref{sec:basics} derive bounds on certain parameters of locally $n \times n$ grid graphs. We look at maximal cliques of locally $n \times n$ grid graphs in Section \ref{sec:cliques}. Finally, in Section \ref{sec:lowerbound} we restrict further to the case where all $\mu$-graphs have order at least $2(n-1)$, and prove Theorem \ref{maintheorem:order-diameter} (\ref{>2(n-1)}) and \ref{maintheorem:order-diameter} (\ref{=2(n-1)}). We apply some of these results to the case where $n = 5$ and prove Theorem \ref{maintheorem:5by5}.

\subsection*{Acknowledgement}

The first and second authors acknowledge the hospitality of the Centre for the Mathematics of Symmetry and Computation of UWA, where this research was carried out. The first author was supported by a Post-doctoral Research Award (FRASDP) of the University of the Philippines. The second author was supported by NSFC (12271524)and NSF of Jiangxi (20224ACB201002, 20212BAB201010). The third author was supported by Australian Research Council grant DP130100106. The authors are grateful to Gordon Royle for pointing out the examples in Construction \ref{example:drg}, to Jonathan Hall for generously sharing his paper \cite{3byq}, and to Dima Pasechnik for drawing our attention to the examples arising from extended generalised quadrangles, and in particular pointing out a gap in an earlier version of this paper. Finally we are grateful to Aart Blokhuis and Andries Brouwer, whose paper \cite{4by4} is the basis of this work and the source of many hours of mathematical joy.

%% file: locallygrid-sec2-prelims-sec3-examples.tex
\section{Preliminaries} \label{sec:prelims}

Let $\Gamma$ be a graph. The \emph{order} $|\V(\Gamma)|$ of $\Gamma$ is the cardinality of $\V(\Gamma)$. For any $x, y \in \V(\Gamma)$, the \emph{distance} $d_\Gamma(x,y)$ in $\Gamma$ of $x$ and $y$ is the length of the shortest path in $\Gamma$ between $x$ and $y$. The \emph{diameter} $\diam(\Gamma)$ of $\Gamma$ is the maximum possible distance between two vertices of $\Gamma$.

Throughout we use the following notation: For $0 \leq i \leq \diam(\Gamma) = D$ and $x \in \V(\Gamma)$ we write $\Gamma_i(x) = \{ y \ : \ d_\Gamma(x,y) = i \}$; we often write $\Gamma(x) = \Gamma_1(x)$. For $y \in \Gamma_i(x)$, \\
	\begin{minipage}[c]{0.48\textwidth}
	\begin{align*}
	k_i(x) &:= |\Gamma_i(x)| \\
	a_i(x,y) &:= |\Gamma_i(x) \cap \Gamma_1(y)| \\
	b_i(x,y) &:= |\Gamma_{i+1}(x) \cap \Gamma_1(y)|, \ 0 \leq i \leq D - 1 \\
	c_i(x,y) &:= |\Gamma_{i-1}(x) \cap \Gamma_1(y)|, \ 1 \leq i \leq D
	\end{align*}
	\end{minipage}
	\hfill
	\begin{minipage}[c]{0.48\textwidth}
	\begin{center}
	\begin{pspicture}(-2.5,-2)(2.5,1)
	\pnode(-2.5,0){C} \pnode(2.5,0){B} \pnode(0,-0.37){K2}
	\rput(0,0){\ovalnode{K}{$k_i(x)$}}
	\ncline[arrows=<-]{C}{K} \ncline[arrows=->]{K}{B} \nccircle[angleA=180,nodesepA=3pt,arrows=->]{K2}{0.25cm}
	\rput(-1.5,0.25){\small $c_i(x,y)$} \rput(1.5,0.25){\small $b_i(x,y)$} \rput(0,-1.2){\small $a_i(x,y)$}
	\end{pspicture}
	\end{center}
	\end{minipage}

In particular, if $\Gamma$ is locally $n \times n$ grid, then $k_1(x) = |K_n \square K_n| = n^2$ for each $x$, and since each vertex in $K_n \square K_n$ has $2(n-1)$ neighbours we have $a_1(x,y) = 2(n-1)$ for each $x \in \V(\Gamma)$ and $y \in \Gamma(x)$. Thus
	\[ b_1(x,y) = k_1(x) - a_1(x,y) - 1 = (n-1)^2. \]
If $d_\Gamma(x,y) = 2$ we usually write $\mu(x,y) = \Gamma(x) \cap \Gamma(y)$ for the $\mu$-graph, and so $c_2(x,y) = |\Gamma(x) \cap \Gamma(y)| = |\mu(x,y)|$. 

In general the parameters $k_i$, $a_i$, $b_i$, and $c_i$ may be non-constant: $k_i(x)$ may depend on $x$, and $a_i(x,y)$, $b_i(x,y)$, and $c_i(x,y)$ may depend on both $x$ and $y$. When they are independent of $x$ or $y$ we sometimes omit the $x$ or $y$. So for example, if $\Gamma$ is locally $n \times n$ grid, we often write $k_1 = n^2$, $a_1 = 2(n-1)$, and $b_1 = (n-1)^2$.

If $b_i$ and $c_i$ are independent of $x$ and $y$ for all $i \in \{0, \ldots, \diam(\Gamma)\}$, then $\Gamma$ is \emph{distance-regular} with \emph{intersection array} $(b_0, b_1, \ldots, b_{D-1}; c_1, c_2, \ldots, c_D)$. In this case the parameters $k_i$ and $a_i$ are also independent of $x$ and $y$, and are determined by the intersection array.

For any $S \subseteq \V(\Gamma)$, we denote by $[S]$ the induced subgraph of $\V(\Gamma)$ on $S$. For any $i \in \{2, \ldots, \diam(\Gamma)\}$, denote by $\Gamma_i$ the set of all pairs of vertices $(x,y)$ such that $d_\Gamma(x,y) = i$. For any $x \in \V(\Gamma)$ define the \emph{eccentricity} $\epsilon(x)$ of $x$ as
	\begin{equation} \label{eq:epsilon}
	\epsilon(x) := \max\{ i \ : \ \Gamma_i(x) \neq \varnothing \}.
	\end{equation}
Clearly $\epsilon(x) \leq \diam(\Gamma)$ for any vertex $x$.

\subsection{Extended generalised quadrangles} \label{ss:egq}

Locally grid graphs arise as point graphs of incidence structures called extended generalised quadrangles. A \emph{generalised quadrangle} (GQ) is an incidence structure of \emph{points} and \emph{lines}, with the following properties: any two distinct points are in at most one common line; if $\ell$ is a line and $P$ a point not on $\ell$ then there is a unique point $Q$ on $\ell$ such that $P$ and $Q$ are collinear; and every object is incident with at least two others. The \emph{point graph} of a GQ is the graph whose vertices are the points of the GQ and whose edges are the pairs of distinct collinear points. For example, an $m \times n$ array of points, with the lines being the rows and columns of the array, is a GQ whose point graph is the grid graph $K_m \square K_n$. Hence any grid graph can be constructed as the point graph of a GQ (but the converse is not true in general). An \emph{extended generalised quadrangle} (EGQ) is an incidence structure of points and lines such that for each point $x$, the \emph{residue} $P_x$ consisting of all points collinear with $x$ and all lines containing $x$ is a GQ under the same incidence relation. The point graph of an EGQ is defined similarly as for a GQ, that is, the vertices are the points of the EGQ, and the lines of the EGQ correspond to maximal cliques in the graph. The residue of a point $x$ in an EGQ corresponds to the induced subgraph on the set of all neigbhours of $x$ in the point graph. If a graph is locally grid then the induced subgraph on the set of all neighbours of a vertex is a grid, which, as already explained, corresponds to a GQ. It follows that a locally grid graph is the point graph of an EGQ.

In the case where $m = n$, the GQ whose lines are the rows and columns of an $n \times n$ array has the additional property that any point lies in two lines and any line contains $n$ points; such a GQ is said to have order $(n-1,1)$. (In general a GQ has order $(s,t)$ if each point lies on $t+1$ lines and each line contains $s+1$ points.) The result \cite[Theorem 2.9]{egq2} gives basic properties of EGQs where each residue $P_x$ is a GQ of order $(s_x,t_x)$ for some integers $s_x$ and $t_x$; \cite[Theorem 2.9 (a)]{egq2} states that in this case the integers $s_x$ and $t_x$ are independent of $x$, that is, there exist integers $s$ and $t$ such that each residue is a GQ of order $(s,t)$. This implies that, if a connected graph $\Gamma$ is locally grid such that for each vertex $x$ the subgraph $[\Gamma(x)]$ is an $n_x \times n_x$ grid for some integer $n_x$, then $\Gamma$ is locally $n \times n$ grid for some integer $n$ that is independent of $x$. The next result lists other elementary properties of locally $n \times n$ grid graphs; parts \eqref{basic-maxcliques}--\eqref{basic-cliques-meet} follow immediately from the other statements in \cite[Theorem 2.9]{egq2}.


\begin{lemma} \label{lemma:basic-nxn}
Let $\Gamma$ be connected and locally $n \times n$ grid. Then:
	\begin{enumerate}[(1)]
	\item \label{basic-maxcliques} A maximal clique in $\Gamma$ has size $n+1$, and each vertex is in $2n$ maximal cliques.
	\item \label{basic-edges} Each edge $\{x,y\}$ is in $2$ maximal cliques and in $2(n-1)$ triangles. Moreover, $[\Gamma(x) \cap \Gamma(y)] \cong 2\,K_{n-1}$.
	\item \label{basic-cliques-meet} Two distinct maximal cliques in $\Gamma$ have $0$ or $2$ vertices in common.
	\item \label{basic-triangles} The number of maximal cliques is $|\V(\Gamma)| \cdot 2n/(n+1)$ and the number of triangles is $|\V(\Gamma)| \cdot n^2(n-1)/3$. Hence $n+1$ divides $2|\V(\Gamma)|$, and if $n \equiv 2\pmod{3}$ then $3$ divides $|\V(\Gamma)|$.
	\item \label{basic-mu} Each $\mu$-graph is a union of $\ell$ cycles, say of lengths $2m_1, \ldots, 2m_\ell$, where each $m_i \geq 2$ and $\sum_{i=1}^\ell m_i \leq n$. No two edges of $\mu(x,y)$ lie in the same $n$-clique in $[\Gamma(x)]$ or $[\Gamma(y)]$.
	\end{enumerate}
\end{lemma}

\begin{proof}
Statements \eqref{basic-maxcliques}--\eqref{basic-cliques-meet} follow immediately from \cite[Theorem 2.9]{egq2}.

By statement (\ref{basic-maxcliques}) each vertex is in $2n$ maximal cliques, and each maximal clique contains $n+1$ vertices. Hence there are $|\V(\Gamma)| \cdot 2n/(n+1)$ maximal cliques. Each vertex is in $n^2$ edges, and by statement (\ref{basic-edges}) each edge is in $2(n-1)$ triangles. Each vertex is contained in two edges in the same triangle, and each triangle has three edges. Therefore the number of triangles is $|\V(\Gamma)| \cdot n^2 \cdot 2(n-1)/6$, and statement (\ref{basic-triangles}) follows.

Let $(x,y) \in \Gamma_2$, and let $z \in \Gamma(x) \cap \Gamma(y)$. Then $x$ and $y$ are vertices in $[\Gamma(z)] \cong K_n \square K_n$, and thus $x$ and $y$ have two common neighbours $u$ and $v$ in $[\Gamma(z)]$. The vertices $u$ and $v$ are non-adjacent in $\Gamma$, and are neighbours of $z$ in $\mu(x,y)$. Hence $\mu(x,y)$ is not a complete graph and has valency $2$, which implies that it is a union of cycles, each of length at least $4$. Now $\mu(x,y)$ is a subgraph of $[\Gamma(x)] \cong K_n \square K_n$; since $\mu(x,y)$ has no triangles, no two of its edges can belong to the same clique of $[\Gamma(x)]$. Thus a connected component of $\mu(x,y)$ has the form given in Figure \ref{figure:mu-component}, and must have even length. Each connected component with length, say, $2m_i$, determines $m_i$ horizontal and $m_i$ vertical cliques, and $m_i \geq 2$. It follows that if $\ell$ is the number of connected components of $\mu(x,y)$, then $\sum_{i=1}^\ell m_i \leq n$. This proves statement (\ref{basic-mu}).
\end{proof}

\begin{center}
\begin{figure}
\begin{pspicture}(-1,-1)(1,1)
\multido{\n=-1+0.4}{5}{\psline[linecolor=lightgray](\n,1)(\n,-1)}
\multido{\n=-0.6+0.4}{5}{\psline[linecolor=lightgray](-1,\n)(1,\n)}
\scalebox{0.8}{\rput(-0.25,0.25){\mucomponent}}
\rput(2.5,0){\small $[\Gamma(x)] \cong K_n \square K_n$}
\end{pspicture}
\caption{Connected component of $\mu(x,y)$} \label{figure:mu-component}
\end{figure}
\end{center}

\section{A family of examples} \label{sec:family}

\begin{construction} \label{example:drg} \cite[Construction 4.1]{DRcovers}
Let $n$ be a power of an odd prime, and let $q = n^2$ and $r = (n+1)/2$. Let $V$ be a vector space of dimension $2$ over the finite field $\F_q$ of order $q$, let $V^*$ be the set of all nonzero vectors, let $B$ be a nondegenerate symplectic form on $V$, and let $R$ be the subgroup of index $r$ in the multiplicative group $\F^*_q$ of $\F_q$. The graph $\Gamma^{(n)}$ has vertex set
	\[ \V\big(\Gamma^{(n)}\big) = \{ Ru \ : \ u \in V^* \} \]
and edge set
	\[ \E\big(\Gamma^{(n)}\big) = \big\{ \{Ru,Rv\} \ : \ B(u,v) \in R \big\}. \]
\end{construction}

By \cite{DRcovers} the graph $\Gamma^{(n)}$ has diameter $3$, and is a distance-regular antipodal cover of $K_{q+1}$ with antipodal blocks of size $r$ and $c_2 = 2(n-1)$. Its intersection array is $\big(q,\, (r-1)c_2,\, 1;\, 1,\, c_2,\, q\big)$. In particular, the graph $\Gamma^{(3)}$ is isomorphic to the Johnson graph $J(6,3)$.

Our aim in this section is to prove Theorem \ref{maintheorem:family}. It will follow from Proposition \ref{proposition:family}.

\begin{proposition} \label{proposition:family}
Let $n$, $q$, $r$, and $\Gamma^{(n)}$ be as in Construction \ref{example:drg}. Then the following hold:
	\begin{enumerate}[(1)]
	\item \label{trans} The graph $\Gamma^{(n)}$ is vertex-transitive and arc-transitive.
	\item \label{locgrid} The graph $\Gamma^{(n)}$ is locally $n \times n$ grid.
	\item \label{cycles-divisors} For each $\mu$-graph of $\Gamma^{(n)}$ there is an odd divisor $d$ of $n-1$ such that the $\mu$-graph is a union of $d$ cycles of length $2(n-1)/d$. Conversely, for each odd divisor $d$ of $n-1$, there is a $\mu$-graph of $\Gamma^{(n)}$ which is a union of $d$ cycles of length $2(n-1)/d$.
	\item \label{cycles-lowerbound} For each $N > 0$ there exists $n \geq N$ such that the $\mu$-graphs of $\Gamma^{(n)}$ are unions of more than $\log(N)$ cycles.
	\end{enumerate}
\end{proposition}

The proof of Proposition \ref{proposition:family} is given at the end of the section, and relies on several intermediate results.

Let $\omega$ be a primitive element of $\F_q$, so that $\omega^{2r}$ is a primitive element of $\F_n$, where $\F_n = \F_{\sqrt{q}} = \left\langle \omega^{2r} \right\rangle \cup \{0\}$ is the subfield of $\F_q$ of index $2$. Then $R = \langle \omega^r \rangle = \F^*_n \,\dot{\cup}\, \F^*_n \omega^r$. The set $\{1, \omega^r\}$ is a basis for $\F_q$ as a vector space over $\F_n$, so $\F_q = \F_n + \F_n \omega^r$ and each $\alpha \in \F_q$ can be written uniquely as
	\begin{equation} \label{eq:E+O}
	\alpha = \alpha_{ev} + \alpha_{odd}, \quad \text{for } \alpha_{ev} \in \F_n \text{ and } \alpha_{odd} \in \F_n\omega^r.
	\end{equation}
Observe that $-1 = \omega^{r(n-1)}$, so $-1 \in \F_n$ (since $n$ is odd) and in particular $-1 \in R$. Also note that $\alpha_{ev}\alpha_{odd}^{-1} \in \F^*_n\omega^r$.

In what follows $\{e,f\}$ is a symplectic basis for $V$ with respect to the form $B$, that is, $e$ and $f$ are nonzero vectors satisfying $B(e,e) = B(f,f) = 0$ and $B(e,f) = -B(f,e) = 1$. Note that $B(u,u) = 0$ for all $u \in V$.

Since $\Gamma^{(n)}$ is an antipodal distance-regular graph of diameter $3$, the antipodal block containing any vertex  $u$ is $\{u\} \cup \Gamma^{(n)}_3(u)$.

\begin{lemma} \label{lemma:antipodes}
Let $\Gamma^{(n)}$ be as in Construction \ref{example:drg}. For any $u \in V^*$, the antipodal block containing $Ru$ is $\{ R'u \ : \ R' \text{ an $R$-coset in $\F_q$} \}$ and this block is $\{Ru\} \cup \Gamma^{(n)}_3(Ru)$.
\end{lemma}

\begin{proof}
Let $R' \neq R$ be an $R$-coset in $\F_q$. Then $R' = R\gamma$ for some $\gamma \notin R$, and $R'u = R(\gamma u)$. Now $B(u,\gamma u) = \gamma B(u,u) = 0 \notin R$, so $R'u \notin \Gamma^{(n)}(Ru)$. Let $Rv \in \Gamma^{(n)}(Ru)$. Then $B(u,v) \in R$ so that $B(\gamma u,v) = \gamma B(u,v) \notin R$, and hence $Rv \notin \Gamma^{(n)}(R(\gamma u)) = \Gamma^{(n)}(R'u)$. Thus $\Gamma^{(n)}(Ru) \cap \Gamma^{(n)}(R'u) = \varnothing$, so that $R'u \notin \Gamma^{(n)}_2(Ru)$. Since $\diam\big(\Gamma^{(n)}\big) = 3$ it follows that $R'u \in \Gamma^{(n)}_3(Ru)$. Therefore $Ru$ and $R'u$ are at maximum distance in $\Gamma^{(n)}$. As mentioned above $\Gamma^{(n)}$ is antipodal and its antipodal blocks have size $r$; since $R$ has index $r$ in $\F^*_q$ the result follows. 
\end{proof}

The action on vectors of the isometry group $\Sp_2(q)$ of $B$ induces an action on $\V\big(\Gamma^{(n)}\big)$ which preserves $\E\big(\Gamma^{(n)}\big)$. This together with the subgroup of scalars isomorphic to $R$ generates $G := R \circ \Sp_2(q)$; again the $G$-action on vectors induces an action on $\V\big(\Gamma^{(n)}\big)$ whose kernel is $R$. (It is convenient to work with this unfaithful action rather than the induced group $\textnormal{PSp}_2(q)$.) We represent vectors in $V$ a row vectors, so $\alpha e + \beta f$ is represented as $(\alpha,\beta)$ and then $G$ acts by matrix multiplication.

\begin{lemma} \label{lemma:family1}
Let $n$, $q$, $r$, and $\Gamma^{(n)}$ be as in Construction \ref{example:drg}. Set $x := Re \in \V\big(\Gamma^{(n)}\big)$ and let $G = R \circ \Sp_2(q)$.
	\begin{enumerate}[(1)]
	\item \label{Gam(x)} $\Gamma^{(n)}(x) = \{ R(\alpha e + f) \ : \alpha \in \F_q \}$, and the stabiliser $G_x$ of $x$ is transitive on $\Gamma^{(n)}(x)$.
	\item \label{Gam(x)-adj} Two distinct vertices $R(\alpha e + f), R(\alpha'e + f) \in \Gamma^{(n)}(x)$ are adjacent in $\Gamma^{(n)}$ if and only if either $\alpha_{ev} = \alpha'_{ev}$ or $\alpha_{odd} = \alpha'_{odd}$ (but not both), where $\alpha_{ev}$, $\alpha'_{ev}$, $\alpha_{odd}$, and $\alpha'_{odd}$ are as in equation (\ref{eq:E+O}). The maximal cliques in $\Gamma^{(n)}(x)$ which contain $R(\alpha e + f)$ are $\big\{ R((\alpha_{ev} + \gamma)e + f) \ : \ \gamma \in \F_n\omega^r \big\}$ and $\big\{ R((\gamma + \alpha_{odd})e + f) \ : \ \gamma \in \F_n \big\}$.
	\item \label{Gam2(x)} $\Gamma^{(n)}_2(x) = \big\{ R(\alpha e + \beta f) \ : \ \alpha \in \F_q,\, \beta \in \F^*_q \setminus R \big\}$. For any $R(\alpha e + \beta f) \in \Gamma^{(n)}_2(x)$, there exists $g \in G_x$ such that $(R(\alpha e + \beta f))^g = R(\beta f)$.
	\end{enumerate}
\end{lemma}

\begin{proof}
Let $\alpha, \beta \in \F_q$. The vertex $R(\alpha e + \beta f)$ is adjacent to $x$ if and only if $\beta = B(e, \alpha e + \beta f) \in R$. Since $0 \notin R$, each such vertex in $\Gamma^{(n)}$ has a unique representative of the form $\alpha e + f$, which proves the first part of statement (\ref{Gam(x)}). For any $\alpha, \alpha' \in \F_q$ the element
	\[ \left(\begin{array}{cc} 1 & 0 \\ \alpha'-\alpha & 1 \end{array}\right) \]
of $G$ fixes $x$ and sends the vertex $R(\alpha e + f)$ to $R(\alpha' e + f)$. This completes the proof of statement (\ref{Gam(x)}).

Let $y = R(\alpha e + f)$ and $z = R(\alpha' e + f)$. Then $y \sim_{\Gamma^{(n)}} z$ if and only if
	\[ \alpha - \alpha' = B(\alpha e + f, \, \alpha'e + f) \in R = \F^*_n \,\dot{\cup}\, \F^*_n\omega^r. \]
Using the representation in equation (\ref{eq:E+O}), $\alpha - \alpha' = (\alpha_{ev} - \alpha'_{ev}) + (\alpha_{odd} - \alpha'_{odd})$. Both $\F_n$ and $\F_n\omega^r$ are closed under addition, so $\alpha_{ev} - \alpha'_{ev} \in \F_n$ and $\alpha_{odd} - \alpha'_{odd} \in \F_n\omega^r$. Thus $\alpha - \alpha' \in \F^*_n$ if and only if $\alpha_{odd} - \alpha'_{odd} \in \F_n$, or equivalently $\alpha_{odd} - \alpha'_{odd} \in \F_n \cap \F_n\omega^r = \{0\}$. Similarly $\alpha - \alpha' \in \F^*_n\omega^r$ if and only if $\alpha_{ev} - \alpha'_{ev} \in \F_n\omega^r$, that is, $\alpha_{ev} - \alpha'_{ev} \in \F_n\omega^r \cap \F_n = \{0\}$. Hence $y \sim_{\Gamma^{(n)}} z$ if and only if either $\alpha_{ev} = \alpha'_{ev}$ or $\alpha_{odd} = \alpha'_{odd}$, but not both (since $y \neq z$). This proves the first part of statement (\ref{Gam(x)-adj}). The second part follows immediately.

The vertex $R(\alpha e + \beta f) \in \Gamma^{(n)}_2(x)$ if and only if $\beta \neq 0$ (for otherwise $R(\alpha e + \beta f) \in \Gamma^{(n)}_3(x)$ by Lemma \ref{lemma:antipodes}) and $\beta \notin R$ (else $R(\alpha e + \beta f) \in \Gamma^{(n)}(x)$ by the above). Hence we obtain the first part of statement (\ref{Gam2(x)}). For any $\alpha \in \F_q$ and $\beta \in \F^*_q \setminus R$ the stabiliser $G_x$ contains the element
	\[ \left(\begin{array}{cc} 1 & 0 \\ -\alpha\beta^{-1} & 1 \end{array}\right), \]
and this sends $R(\alpha e + \beta f)$ to $R(\beta f)$. This completes the proof of statement (\ref{Gam2(x)}).
\end{proof}

By Lemma \ref{lemma:family1} (\ref{Gam2(x)}) each $G_x$-orbit in $\Gamma^{(n)}_2(x)$ contains a vertex $R\beta f$ for some $\beta \in \F^*_q \setminus R$. Hence to determine the structure of the $\mu$-graphs $\mu(x,y)$ for any $y$ we may assume that $y \in \Gamma^{(n)}_2(x)$ is $R\beta f$. This is what we do in the next result.

We denote the multiplicative order of $\alpha \in \F^*_q$ by $|\alpha|$.

\begin{lemma} \label{lemma:family2}
Let $n$, $q$, $r$, and $\Gamma^{(n)}$ be as in Construction \ref{example:drg}. Set $x := Re \in \V\big(\Gamma^{(n)}\big)$ and $y = R(\beta^{-1}f)$ where $\beta \in \F^*_q \setminus R$, and let $G = R \circ \Sp_2(q)$.
	\begin{enumerate}[(1)]
	\item \label{mu(x,y)-adj} $y \in \Gamma^{(n)}_2(x)$ and $\Gamma^{(n)}(x) \cap \Gamma^{(n)}(y) = \{ R(\alpha e + f) \ : \ \alpha \in R\beta \}$ of size $2(n-1)$. Two distinct vertices $R(\alpha e + f), R(\alpha'e + f) \in \Gamma^{(n)}(x) \cap \Gamma^{(n)}(y)$ are adjacent in $\Gamma^{(n)}$ if and only if $\alpha' = \alpha\big(\beta_{ev}\beta_{odd}^{-1}\big)^{\pm 1}$.
	\item \label{mu(x,y)-cycles} The $\mu$-graph $\mu(x,y)$ is a union of $d := 2(n-1)/\left|\beta_{ev}\beta_{odd}^{-1}\right|$ cycles of length $2(n-1)/d$.
	\end{enumerate}
\end{lemma}

\begin{proof}
Since $\beta \notin R \cup \{0\}$ neither is $\beta^{-1}$, so $y \in \Gamma^{(n)}_2(x)$ by Lemma \ref{lemma:family1} (\ref{Gam2(x)}). For any $\alpha', \beta' \in \F_q$, we have $R(\alpha'e + \beta'f) \in \Gamma^{(n)}(y)$ if and only if $-\beta^{-1}\alpha' = B(\beta^{-1}f, \, \alpha'e + \beta'f) \in R$, or equivalently $\alpha' \in R\beta$ (since $-1 \in R$). Thus $\Gamma^{(n)}(y) = \big\{ R(\alpha'e + \beta'f) \ : \ \alpha' \in R\beta, \, \beta' \in \F_q \big\}$, and we conclude that
	\[ \Gamma^{(n)}(x) \cap \Gamma^{(n)}(y) = \big\{ R(\alpha' e + \beta' f) \ : \ \alpha' \in R\beta, \, \beta' \in R \big\},
	= \big\{ R(\alpha' e + f) \ : \ \alpha' \in R\beta \big\}. \]
a set of size $|R\beta| = 2(n-1)$. This proves the first part of statement (\ref{mu(x,y)-adj}).

Next let $w_1 = R(\alpha_1 e + f)$ and $w_2 = R(\alpha_2 e + f)$ be distinct vertices in $\Gamma^{(n)}(x) \cap \Gamma^{(n)}(y)$. Then $\alpha_1$ and $\alpha_2$ are distinct elements of $R\beta$, and for each $i \in \{1,2\}$ we can write $\alpha_i = \rho_i\beta$ for some $\rho_i \in R$. Hence $\alpha_i = \rho_i(\beta_{ev} + \beta_{odd})$ for $i = 1, 2$. If both $\rho_1, \rho_2 \in \F^*_n$ then $(\alpha_i)_{ev} = \rho_i\beta_{ev}$ and $(\alpha_i)_{odd} = \rho_i\beta_{odd}$, and since $\alpha_1 \neq \alpha_2$ either $(\alpha_1)_{ev} \neq (\alpha_2)_{ev}$ or $(\alpha_1)_{odd} \neq (\alpha_2)_{odd}$. So $\rho_1 \neq \rho_2$, and both $(\alpha_1)_{ev} \neq (\alpha_2)_{ev}$ and $(\alpha_1)_{odd} \neq (\alpha_2)_{odd}$ hold. Thus $w_1 \nsim_\Gamma^{(n)} w_2$ by Lemma \ref{lemma:family1} (\ref{Gam(x)-adj}), and we can also deduce by a similar argument that $w_1 \nsim w_2$ whenever both $\rho_1, \rho_2 \in \F^*_n \omega^r$. Let us therefore assume that $\rho_1$ and $\rho_2$ belong in different $\F^*_n$-cosets in $R$; without loss of generality suppose that $\rho_1 \in \F^*_n$ and $\rho_2 \in \F^*_n\omega^r$. Then $(\alpha_1)_{ev} = \rho_1\beta_{ev}$, $(\alpha_1)_{odd} = \rho_1\beta_{odd}$, $(\alpha_2)_{ev} = \rho_2\beta_{odd}$, and $(\alpha_2)_{odd} = \rho_2\beta_{ev}$. By Lemma \ref{lemma:family1} (\ref{Gam(x)-adj}) the vertices $w_1$ and $w_2$ are adjacent if and only if either $\rho_1\beta_{ev} = \rho_2\beta_{odd}$ or $\rho_1\beta_{odd} = \rho_2\beta_{ev}$ (but not both), which is equivalent to $\rho_2 = \rho_1\left(\beta_{ev}\beta_{odd}^{-1}\right)^{\pm 1}$. So $w_1 \sim_{\Gamma^{(n)}} w_2$ if and only if
	\[ \alpha_2 = \rho_1\left(\beta_{ev}\beta_{odd}^{-1}\right)^{\pm 1}\beta = \alpha_1\big(\beta_{ev}\beta_{odd}^{-1}\big)^{\pm 1}. \]
This completes the proof of statement (\ref{mu(x,y)-adj}). (Recall that $\beta_{ev}\beta^{-1}_{odd} \in \F^*_n \omega^r \subseteq R$ by equation (\ref{eq:E+O}).)

Set $\gamma = \beta_{ev}\beta_{odd}^{-1}$. It follows from the above that
	\[ R(\alpha e + f) \sim_{\Gamma^{(n)}} R(\alpha\gamma e + f) \sim_{\Gamma^{(n)}} R(\alpha\gamma^2 e + f) \sim_{\Gamma^{(n)}} \ldots, \]
that is, each connected component of $\mu(x,y)$ has vertex set $\big\{ R(\alpha' e + f) \ : \ \alpha' \in \alpha \big\langle \beta_{ev}\beta_{odd}^{-1} \big\rangle \big\}$ for some $\alpha \in \F_q$. Note that $R(\alpha\gamma^k e + f) = R(\alpha e + f)$ by the uniqueness of the coset representative of the form $\alpha e + f$. Thus the length of each component is $|\gamma|$, and the number of components is $d = 2(n-1)/|\gamma|$. This proves statement (\ref{mu(x,y)-cycles}).
\end{proof}

\begin{proof}[Proof of Proposition \ref{proposition:family}]
Let $G = R \circ \Sp_2(q)$. Then the $G$-action on $V^*$ induces an action on $\V\big(\Gamma^{(n)}\big)$ which preserves $\E\big(\Gamma^{(n)}\big)$, and the kernel of this action is $R$. Since $\Sp_2(q)$ acts transitively on $V^*$, the group $G$ is transitive on $\V\big(\Gamma^{(n)}\big)$. By Lemma \ref{lemma:family1} (\ref{Gam(x)}) the stabiliser in $G$ of the vertex $x = Re$ is transitive on $\Gamma^{(n)}(x)$; since $G$ is vertex-transitive, it follows that $G$ is also arc-transitive on $\Gamma^{(n)}$. This proves statement (\ref{trans}).

It is easy to see from Lemma \ref{lemma:family1} (\ref{Gam(x)-adj}) that $[\Gamma^{(n)}(x)] \cong K_n \square K_n$, so by vertex-transitivity $\Gamma^{(n)}$ is locally $n \times n$ grid. Hence statement (\ref{locgrid}) holds.

To prove statement \ref{cycles-divisors} first let $x', y' \in \V\big(\Gamma^{(n)}\big)$ with $d_\Gamma^{(n)}(x',y') = 2$. By vertex-transitivity and Lemma \ref{lemma:family1} (\ref{Gam2(x)}) there exist $g \in G$ and $h \in G_x$ such that $(x')^g = x = Re$ and $(y')^{gh} = R\beta^{-1}f =: y$ for some $\beta \in \F^*_q \setminus R$. That is, $(x',y')^{gh} = (x,y)$, so that $\mu(x',y') \cong \mu(x,y)$. By Lemma \ref{lemma:family2} (\ref{mu(x,y)-cycles}), the graph $\mu(x,y)$ is a union of $d = 2(n-1)/\big|\beta_{ev}\beta_{odd}^{-1}\big|$ cycles of length $\big|\beta_{ev}\beta_{odd}^{-1}\big|$. Since $\beta_{ev}\beta_{odd}^{-1} \in \F^*_n\omega^r$, we have $\beta_{ev}\beta_{odd}^{-1} = \omega^{ri}$ for some odd $i$. Thus $\big|\beta_{ev}\beta_{odd}^{-1}\big| = (q-1)/\gcd(ri,q-1) = 2(n-1)/\gcd(i,2(n-1))$, implying that $d = \gcd(i,2(n-1))$. Further, since $i$ is odd, $d = \gcd(i,n-1)$ is odd. Thus $d$ is an odd divisor of $n-1$. This proves the first part of statement (\ref{cycles-divisors}).

For the converse, let $d$ be an odd divisor of $n-1$. Take $x = Re$ and $y = R\beta^{-1}f$, where $\beta = 1 + \omega^{-rd}$. Note that $\omega^{-rd} = \omega^{r(q-1-d)}$; since both $q$ and $d$ are odd, so is $q-1-d$, and thus $\omega^{-rd} \in \F^*_n\omega^r$. Hence $\beta_{ev} = 1$ and $\beta_{odd} = \omega^{-rd}$. Since for $\gamma = 0$ we have $\gamma_{ev} = \gamma_{odd} = 0$, it follows from the uniqueness of the expression (\ref{eq:E+O}) for $\beta$ that $\beta \neq 0$. Also $\beta \notin \F^*_n \cup \F^*_n\omega^r = R$ since $\beta_{ev}$ and $\beta_{odd}$ are both nonzero. Therefore $\beta \in \F^*_q \setminus R$, so $y \in \Gamma^{(n)}_2(x)$ by Lemma \ref{lemma:family1} (\ref{Gam2(x)}). Now $\big|\beta_{ev}\beta_{odd}^{-1}\big| = \big|\omega^{rd}\big| = 2(n-1)/d$, so by Lemma \ref{lemma:family2} (\ref{mu(x,y)-cycles}) the graph $\mu(x,y)$ is a union of $d$ cycles of length $2(n-1)/d$. This completes the proof of statement (\ref{cycles-divisors}).

It follows from Proposition \ref{proposition:family} (\ref{cycles-divisors}) that there is no absolute upper bound on the number of cycles in a $\mu$-graph in Construction \ref{example:drg}. For if $n = p^m$ for some odd prime $p$, and $m \geq 3$, then $p^m - 1$ has a prime divisor $d$ that does not divide $p^i - 1$ for $i < m$ by \cite{ppd-zsig} (see also \cite[Theorem 2.1]{ppd}), and such a prime is at least $m+1 > \log_p(n)$. This proves statement (\ref{cycles-lowerbound}).
\end{proof}

\begin{proof}[Proof of Theorem \ref{maintheorem:family}]
As mentioned before Proposition \ref{proposition:family}, $\diam\big(\Gamma^{(n)}\big) = 3$. Also $\Gamma^{(n)}$ is an antipodal cover of $K_{q+1}$ with antipodal blocks of size $r$, so that $|\V\big(\Gamma^{(n)}\big)| = (q+1)r = (n^2+1)(n+1)/2$. It is distance-regular with parameter $c_2 = 2(n-1)$, so all of its $\mu$-graphs have order $2(n-1)$, and its intersection array is $\big(q,\, (r-1)c_2,\, 1;\, 1,\, c_2,\, q\big) = \big(n^2,\, (n-1)^2,\, 1;\, 1,\, 2(n-1),\, n^2\big)$. It is locally $n \times n$ grid by Proposition \ref{proposition:family} (\ref{locgrid}). Thus $\Gamma^{(n)}$ satisfies the conditions of Theorem \ref{maintheorem:order-diameter} (\ref{=2(n-1)}). The last part of Theorem \ref{maintheorem:family} follows immediately from Proposition \ref{proposition:family} (\ref{cycles-divisors}).
\end{proof}

%% file: locallygrid-sec4-basics.tex
\section{Basic properties of locally $n \times n$ grid graphs} \label{sec:basics}

In this section we establish some basic properties of locally $n \times n$ grid graphs and prove the first statement of Theorem \ref{maintheorem:order-diameter}.

The first result follows from Theorems 3.4 and 3.16 in \cite{egq2}. It gives a tight upper bound on the diameter of a locally $n \times n$ grid graph, and characterises the Johnson graphs as locally grid graphs with maximum diameter.

\begin{theorem} \cite[Theorems 3.4 and 3.16(a)]{egq2} \label{theorem:diameter}
If $\Gamma$ is a locally $n \times n$ grid graph then $\diam(\Gamma) \leq n$. Moreover, $\diam(\Gamma) = n$ if and only if $\Gamma$ is the Johnson graph $J(2n,n)$.
\end{theorem}

\begin{remark} \label{remark:distancediagram}
Let $x \in \V(\Gamma)$, with eccentricity $\epsilon(x)$ as in (\ref{eq:epsilon}), and $2 \leq i \leq \epsilon(x)$. Counting in two ways the number of edges between $\Gamma_{i-1}(x)$ and $\Gamma_i(x)$ yields the equality
	\begin{equation} \label{eq:b=c}
	\sum_{y \in \Gamma_{i-1}(x)} b_{i-1}(x,y) = \sum_{z \in \Gamma_i(x)} c_i(x,z).
	\end{equation}
By Lemma \ref{lemma:basic-nxn} (\ref{basic-mu}), for any $z \in \Gamma_2(x)$ we have $c_2(x,z) = 2m$ for some $m \in \{2, \ldots, n\}$, where $m$ may depend on $x$ and $z$. For $2 \leq m \leq n$ define
	\begin{equation} \label{eq:k_2,2m}
	k_{2,2m}(x) := \big|\big\{ z \in \Gamma_2(x) \ : \ c_2(x,z) = 2m \big\}\big|,
	\end{equation}
so that
	\begin{equation} \label{eq:k2-sum}
	k_2(x) = \sum_{m=2}^n k_{2,2m}(x).
	\end{equation}
Also $\sum_{z \in \Gamma_2(x)} c_2(x,z) = \sum_{m=2}^n 2m\,k_{2,2m}(x)$, and since $k_1(x) = n^2$ and $b_1(x,y) = (n-1)^2$ for all $y \in \Gamma(x)$, we have $\sum_{y \in \Gamma(x)} b_1(x,y) = n^2(n-1)^2$. Thus, for $i = 2$, equation (\ref{eq:b=c}) becomes
	\begin{equation} \label{eq:sum-k_2,2m}
	n^2(n-1)^2 = \sum_{m=2}^n 2m\,k_{2,2m}(x).
	\end{equation}
\end{remark}

\begin{lemma} \label{lemma:maxcliques}
Assume that $\Gamma$ is locally $n \times n$ grid, and let $(x,y) \in \Gamma_2$. Then $c_2(x,y) = 2m$ for some $m \in \{2, \ldots, n \}$, and the following hold:
	\begin{enumerate}[(1)]
	\item \label{2n} If $c_2(x,y) = 2n$ then $d_\Gamma(x,C) = 1$ for any maximal clique $C$ containing $y$.
	\item \label{2m} If $c_2(x,y) = 2m \leq 2(n-1)$, then, of the $2n$ maximal cliques $C$ containing $y$, $d_\Gamma(x,C) = 1$ for $2m$ cliques and $d_\Gamma(x,C) = 2$ for the remaining $2(n-m)$ cliques.
	\end{enumerate}
\end{lemma}

\begin{proof}
By Lemma \ref{lemma:basic-nxn} (\ref{basic-mu}), $c_2(x,y) = 2m$ for some $m \in \{2, \ldots, n\}$ and no two edges of $\mu(x,y)$ lie in the same $n$-clique of $[\Gamma(y)]$. So in $[\Gamma(y)]$ there are $m$ horizontal $n$-cliques and $m$ vertical $n$-cliques that contain an edge of $\mu(x,y)$, as illustrated on the left in Figure \ref{figure:maxcliques-1}. If $m = n$ then each $n$-clique in $[\Gamma(y)]$ contains an edge of $\mu(x,y)$, so that each $(n+1)$-clique containing $y$ is adjacent to $x$, proving statement (\ref{2n}). If $m < n$ then the remaining $n-m$ horizontal $n$-cliques and $n-m$ vertical $n$-cliques in $[\Gamma(y)]$ do not contain any vertex of $\mu(x,y)$, as illustrated on the right in Figure \ref{figure:maxcliques-1}, but each of these cliques contains at least one vertex that is adjacent to a vertex of $\mu(x,y)$. Hence $d_\Gamma(x,C) = 2$ for these cliques $C$, as required.
\end{proof}

\begin{center}
\begin{figure}
\begin{pspicture}(-6,-2.25)(6,2.5)
\rput(-4.75,0){
	\localgrid
	\rput(-0.75,0.75){\scalebox{0.55}{\mucomponent}}
	\rput(0,1.9){\small $[\Gamma(y)]$}
	\rput(-2.2,0.5){\small $\mu(x,y)$} \psline[linewidth=0.25pt](-1.6,0.5)(-1.35,0.5)
	\pnode(1.65,1.5){A} \pnode(1.65,0){B} \ncbar[nodesep=0pt,angle=0,armA=0.1cm,armB=0.1cm]{A}{B} \rput(3.5,0.75){\parbox{3cm}{\small $m$ horizontal \\ cliques}}
	\pnode(-1.5,-1.65){C} \pnode(0,-1.65){D} \ncbar[nodesep=0pt,angle=-90,armA=0.1cm,armB=0.1cm]{C}{D} \rput(-0.75,-2.1){\small $m$ vertical cliques}	
}
\rput(3.25,0){
	\psline[linecolor=lightgray](-1.5,0)(-1.5,1.5)(0,1.5)
	\pspolygon[linecolor=lightgray,fillstyle=vlines,hatchcolor=lightgray,hatchangle=0](0,1.5)(1.5,1.5)(1.5,-1.5)(0,-1.5)
	\pspolygon[linecolor=lightgray,fillstyle=hlines,hatchcolor=lightgray,hatchangle=0](-1.5,0)(-1.5,-1.5)(1.5,-1.5)(1.5,0)
	\rput(-0.75,0.75){\scalebox{0.55}{\mucomponent}}
	\rput(0,1.9){\small $[\Gamma(y)]$}
	\rput(-2.2,0.5){\small $\mu(x,y)$} \psline[linewidth=0.25pt](-1.6,0.5)(-1.35,0.5)
	\pnode(1.65,0){A} \pnode(1.65,-1.5){B} \ncbar[nodesep=0pt,angle=0,armA=0.1cm,armB=0.1cm]{A}{B} \rput(3.5,-0.75){\parbox{3cm}{\small $n-m$ horizontal \\ cliques}}
	\pnode(0,-1.65){C} \pnode(1.5,-1.65){D} \ncbar[nodesep=0pt,angle=-90,armA=0.1cm,armB=0.1cm]{C}{D} \rput(0.75,-2.1){\small $n-m$ vertical cliques}
}
\end{pspicture}
\caption{Maximal cliques $C$ satisfying $d_\Gamma(x,C) = 1$ (left) and $d_\Gamma(x,C) = 2$ (right)} \label{figure:maxcliques-1}
\end{figure}
\end{center}

Statement (\ref{distance1}) of the next lemma is the third assertion in \cite[Lemma, Section 1]{4by4}. The second part of statement (\ref{distance2}) generalises the first assertion in \cite[Lemma 2, Section 5]{4by4}.

\begin{lemma} \label{lemma:maxcliques2}
Assume that $\Gamma$ is locally $n \times n$ grid. Let $x \in \V(\Gamma)$ and $C$ a maximal clique in $\Gamma$ not containing $x$.
	\begin{enumerate}[(1)]
	\item \label{distance1} If $d_\Gamma(x,C) = 1$ then $|C \cap \Gamma(x)| = 2$ and $|C \cap \Gamma_2(x)| = n-1$.
	\item \label{distance2} If $d_\Gamma(x,C) = 2$ then each $y \in C \cap \Gamma_2(x)$ satisfies $c_2(x,y) \leq 2(n-1)$, and if $c_2(x,y) = 2m$ then $|C \cap \Gamma_2(x)| \geq m+1$.
	\end{enumerate}
\end{lemma}

\begin{proof}
Suppose that $d_\Gamma(x,C) = 1$ and that some vertex $y \in C \cap \Gamma(x)$. Then $x \in \Gamma(y)$, and $C \setminus \{y\}$ is an $n$-clique in $\Gamma(y)$ not containing $x$. We see from the $n \times n$ grid $[\Gamma(y)]$ that $x$ is adjacent to a unique vertex in $C \setminus \{y\}$ and is at distance two from any other vertex in $C \setminus \{y\}$. Therefore $|C \cap \Gamma(x)| = 2$ and $|C \cap \Gamma_2(x)| = n-1$, which proves statement (\ref{distance1}).

Suppose now that $d_\Gamma(x,C) = 2$. Let $y \in C \cap \Gamma_2(x)$ and let $C' = C \setminus \{y\}$. Then $c_2(x,y) = 2m \leq 2(n-1)$ (because otherwise by Lemma \ref{lemma:maxcliques} (\ref{2n}) all cliques containing $y$ are at distance $1$ from $x$, and in particular $d_\Gamma(x,C) = 1$, contradiction), and $C'$ is an $n$-clique in $\Gamma(y)$. Since $C' \subseteq C$ we have $d_\Gamma(x,C') \geq 2$, so $C'$ does not contain an edge of $\mu(x,y)$. In $[\Gamma(y)]$ there are $n-1$ cliques of size $n$ that are disjoint from $C'$, and $n$ cliques of size $n$ that meet $C'$ in a unique vertex. Of the $n$ cliques that intersect $C'$, there are $m$ cliques $C_1, \ldots, C_m$ each of which cointains an edge of $\mu(x,y)$. Any two of the cliques $C_i$ are disjoint, so that $|C' \cap (C_1 \cup \ldots \cup C_m)| = m$ as illustrated in Figure \ref{figure:maxcliques-2}. Thus $|C \cap \Gamma_2(x)| \geq |(C' \cap (C_1 \cup \ldots \cup C_m)) \cup \{y\}| = m+1$, which proves statement (\ref{distance2}).
\end{proof}

\begin{center}
\begin{figure}
\begin{pspicture}(-2.5,-1.8)(2.5,2.5)
\multido{\n=-1.5+0.4}{8}{\psline[linecolor=lightgray](\n,1.5)(\n,-1.5)}
\multido{\n=-1.3+0.4}{8}{\psline[linecolor=lightgray](1.5,\n)(-1.5,\n)}
\psline[linewidth=1.3pt,linecolor=lightgray](-1.5,-0.9)(1.5,-0.9) \multido{\n=-1.1+0.4}{5}{\qdisk(\n,-0.9){2pt}}
\rput(-2,-0.9){\small $C'$}
\rput(-0.3,0.3){\scalebox{0.8}{\mucomponent}}
\rput(0,1.9){\small $[\Gamma(y)]$}
\rput(-2.3,0.25){\small $\mu(x,y)$} \psline[linewidth=0.25pt](-1.7,0.25)(-1.2,0.25)
\pnode(-1.1,-1.5){a} \pnode(0.5,-1.5){b} \ncbar[nodesep=2pt,angle=-90,armA=0.1cm,armB=0.1cm]{a}{b} \rput(-0.3,-1.9){\small $m$}
\end{pspicture}
\caption{Vertices in $\Gamma_2(x) \cap C'$} \label{figure:maxcliques-2}
\end{figure}
\end{center}

The next result generalises \cite[Lemmas 1 (iv) and 2]{4by4}.

\begin{lemma} \label{lemma:parameters}
Assume that $\Gamma$ is locally $n \times n$ grid. Let $x \in \V(\Gamma)$ with eccentricity $\epsilon(x)$ as in (\ref{eq:epsilon}).
	\begin{enumerate}[(1)]
	\item \label{b2} Let $y \in \Gamma_2(x)$. If $c_2(x,y) = 2m$ then $b_2(x,y) \leq (n-m)^2$. In particular, if $c_2(x,y) = 2n$ then $b_2(x,y) = 0$.
	\item \label{b3,c3} Assume that $\epsilon(x) \geq 3$ and let $z \in \Gamma_3(x)$. If $c_2(x,y) \geq 2m$ for all $y \in \Gamma_2(x)$ then $c_3(x,z) \geq (m+1)^2$ and $b_3(x,z) \leq (n-m-1)^2$.
	\item \label{bi,ci} Assume that $\epsilon(x) \geq i \geq 4$ and let $z \in \Gamma_i(x)$. If $c_2(x,y) \geq 2m$ for all $y \in \Gamma_2(x)$ then $c_i(x,z) \geq (m+1)^2$ and $b_i(x,z) \leq (n-m-1)^2$.
	\end{enumerate}
\end{lemma}

\begin{proof}
Let $y \in \Gamma_2(x)$ and suppose that $c_2(x,y) = 2m$. Then $\Gamma_3(x) \cap \Gamma(y)$ is contained in the set of all vertices in $\Gamma(y)$ that are not adjacent to any vertex in $\mu(x,y)$, as illustrated in Figure \ref{figure:b,c}. Hence $\Gamma_3(x) \cap \Gamma(y)$ lies in an $(n-m) \times (n-m)$ subgrid of $[\Gamma(y)]$, so it follows that $b_2(x,y) = |\Gamma_3(x) \cap \Gamma(y)| \leq (n-m)^2$. This proves statement (\ref{b2}).

Let $z \in \Gamma_3(x)$ and assume that $c_2(x,y) \geq 2m$ for all $y \in \Gamma_2(x)$. Then $d_\Gamma(x,C) = 2$ for some $(n+1)$-clique $C$ containing $z$. Hence $|C \cap \Gamma_2(x)| \geq m+1$ by Lemma \ref{lemma:maxcliques2} (\ref{distance2}). Let $C' = C \setminus \{z\}$, which is an $n$-clique in $\Gamma(z)$. Since $z \notin \Gamma_2(x) \cap C$, we also have $|C' \cap \Gamma_2(x)| \geq m+1$. Without loss of generality suppose that $C'$ is a ``horizontal'' $n$-clique. Then any ``vertical'' $n$-clique $C''$ containing a point in $C' \cap \Gamma_2(x)$ also satisfies $|C'' \cap \Gamma_2(x)| \geq m+1$, as illustrated in Figure \ref{figure:b,c}. Therefore $c_3(x,z) = |\Gamma_2(x) \cap \Gamma(z)| \geq (m+1)^2$, and $[\Gamma_2(x) \cap \Gamma(z)]$ contains an $(m+1) \times (m+1)$ subgrid. No vertex in $\Gamma_4(x) \cap \Gamma(z)$ is adjacent to any vertex in $\Gamma_2(x) \cap \Gamma(z)$, hence $\Gamma_4(x) \cap \Gamma(z)$ lies in an $(n-m-1) \times (n-m-1)$ subgrid of $\Gamma(z)$. So $b_3(x,z) = |\Gamma_4(x) \cap \Gamma(z)| \leq (n-m-1)^2$, and statement (\ref{b3,c3}) holds.

Finally, let $z \in \Gamma_i(x)$, $4 \leq i \leq \epsilon(x)$, and suppose that $c_2(x,y) \geq 2m$ for all $y \in \Gamma_2(x)$. Take $w \in \Gamma_{i-3}(x)$ such that $d_\Gamma(w,z) = 3$. Then by part \ref{b3,c3} above $c_3(w,z) \geq (m+1)^2$ and $b_3(w,z) \leq (n-m-1)^2$. Clearly $\Gamma_2(w) \cap \Gamma(z) \subseteq \Gamma_{i-1}(x) \cap \Gamma(z)$, and $\Gamma_4(w) \cap \Gamma(z) \subseteq \Gamma_{i+1}(x) \cap \Gamma(z)$. statement (\ref{bi,ci}) follows.
\end{proof}

\begin{center}
\begin{figure}
\begin{pspicture}(-6,-2.5)(6,2.5)
\rput(-4.5,0){
	\localgrid
	\psellipse(0.75,-0.75)(0.6,0.5)
	\rput(-0.75,0.75){\scalebox{0.55}{\mucomponent}}
	\rput(0,1.9){\small $[\Gamma(y)]$}
	\rput(-2.2,0.5){\small $\mu(x,y)$} \psline[linewidth=0.25pt](-1.6,0.5)(-1.35,0.5)
	\rput(2.7,-0.75){\small $\Gamma_3(x) \cap \Gamma(y)$} \psline[linewidth=0.25pt](0.75,-0.75)(1.6,-0.75)
	\pnode(1.65,1.5){A} \pnode(1.65,0){B} \ncbar[nodesep=0pt,angle=0,armA=0.1cm,armB=0.1cm]{A}{B} \rput(3.2,0.75){\parbox{2.5cm}{\small $m$ horizontal \\ cliques}}
	\pnode(-1.5,-1.65){C} \pnode(0,-1.65){D} \ncbar[nodesep=0pt,angle=-90,armA=0.1cm,armB=0.1cm]{C}{D} \rput(-0.75,-2.1){\small $m$ vertical cliques}
}
\rput(3.5,0){
	\localgrid
	\psellipse(0.75,-0.75)(0.6,0.5)
	\multido{\n=0.15+0.40}{4}{\rput(0,\n){\multido{\n=-1.35+0.40}{4}{\qdisk(\n,0){2pt}}}}
	\pnode(-1.5,1.5){A} \pnode(-1.5,0){B} \ncbar[nodesep=4pt,angle=180,armA=0.1cm,armB=0.1cm]{A}{B} \rput(-2.4,0.8){\small $\geq m+1$}
	\pnode(0,1.5){C} \ncbar[nodesep=4pt,angle=90,armA=0.1cm,armB=0.1cm]{A}{C} \rput(-0.8,2){\small $\geq m+1$}
	\psline[linewidth=0.25pt](-1.35,0.15)(-1.7,-0.25) \rput(-2.4,-0.35){\small $\Gamma_2(x) \ni$}
	\rput(2.7,-0.75){\small $\Gamma_4(x) \cap \Gamma(z)$} \psline[linewidth=0.25pt](0.75,-0.75)(1.6,-0.75)
	\rput(0,-1.9){\small $[\Gamma(z)]$}
}
\end{pspicture}
\caption{$[\Gamma(y)]$ for $y \in \Gamma_2(x)$ and $[\Gamma(z)]$ for $z \in \Gamma_3(x)$} \label{figure:b,c}
\end{figure}
\end{center}

\begin{lemma} \label{lemma:k-bounds}
Assume that $\Gamma$ is locally $n \times n$ grid. Suppose that there exists a constant $m \in \{2, \ldots, n-1\}$ such that $c_2(x,y) \geq 2m$ for all $(x,y) \in \Gamma_2$. Then for any $x \in \V(\Gamma)$,
	\begin{align}
	k_2(x) &\leq \frac{n^2(n-1)^2}{2m}, \label{eq:k2} \\
	k_3(x) &\leq k_2(x) \cdot \frac{(n-m)^2}{(m+1)^2} \leq \frac{n^2(n-1)^2(n-m)^2}{2m(m+1)^2}, \label{eq:k3} \\
	\intertext{and for any $4 \leq i \leq \epsilon(x)$, with $\epsilon$ as in (\ref{eq:epsilon}),}
	k_i(x) &\leq k_{i-1}(x) \cdot \frac{(n-m-1)^2}{(m+1)^2}. \label{eq:ki}
	\end{align}
\end{lemma}

\begin{proof}
It follows from the hypothesis and equation (\ref{eq:b=c}) with $i = 2$ that
	\[ n^2(n-1)^2 = \sum_{y \in \Gamma(x)} b_1(x,y) = \sum_{z \in \Gamma_2(x)} c_2(x,z) \geq k_2(x) \cdot 2m, \]
which then yields (\ref{eq:k2}). By Lemma \ref{lemma:parameters} \eqref{b2} and \eqref{b3,c3} we have $b_2(x,y) \leq (n-m)^2$ and $c_3(x,z) \geq (m+1)^2$ for any $y \in \Gamma_2(x)$ and $z \in \Gamma_3(x)$. So with $i = 3$ in (\ref{eq:b=c}) we obtain
	\[ k_2(x) \cdot (n-m)^2 \geq \sum_{y \in \Gamma_2(x)} b_2(x,y) = \sum_{z \in \Gamma_3(x)} c_3(x,z) \geq k_3(x) \cdot (m+1)^2, \]
which then yields (\ref{eq:k3}). Similarly, if $4 \leq i \leq \epsilon(x)$, then by Lemma \ref{lemma:parameters} \eqref{b3,c3} and \eqref{bi,ci} we have $b_{i-1}(x,y) \leq (n-m-1)^2$ and $c_i(x,z) \geq (m+1)^2$ for any $y \in \Gamma_{i-1}(x)$ and $z \in \Gamma_i(x)$. So (\ref{eq:b=c}) gives us
	\[ k_{i-1}(x) \cdot (n-m-1)^2 \geq \sum_{y \in \Gamma_{i-1}(x)} b_{i-1}(x,y) = \sum_{z \in \Gamma_i(x)} c_i(x,z) \geq k_i(x) \cdot (m+1)^2, \]
and (\ref{eq:ki}) follows.
\end{proof}

The distance diagram of a locally grid graph $\Gamma$ satisfying the hypotheses of Lemma \ref{lemma:k-bounds} is shown in Figure \ref{figure:distancediagram}.

\begin{center}
\begin{figure}[ht]
\begin{pspicture}(-7,-0.7)(7,1.25)
\distancediagram
\rput(-2.6,0.25){\small $\geq 2m$} \rput(0.3,0.25){\small $\leq (n-m)^2$}
\rput(2.45,0.25){\small $\geq (m+1)^2$} \rput(6,0.25){\small $\leq (n-m-1)^2$}
\end{pspicture} \\
\begin{pspicture}(-7,-1.25)(7,0.7)
\distancediagramsmall
\rput(-1.6,0.25){\small $\geq (m+1)^2$}
\rput(2,0.25){\small $\leq (n-m-1)^2$}
\end{pspicture}
\caption{Distance diagram for $\Gamma$ with respect to the vertex $x$, assuming $c_2(x',y') \geq 2m$ for all $\{x',y'\} \in \Gamma_2$}
\label{figure:distancediagram}
\end{figure}
\end{center}

%% file: locallygrid-sec5-cliques.tex
\section{Results on maximal cliques} \label{sec:cliques}

In this section we prove some technical results on maximal cliques of locally $n \times n$ grid graphs.

Lemma \ref{lemma:mu-clique} and Corollary \ref{corollary:mu-clique} generalise the result in the first part of the proof of \cite[Lemma 6]{4by4}, and use very similar arguments. If $C$ is a maximal clique in $\Gamma$ and $x$ is a vertex not contained in $C$, it follows from Lemma \ref{lemma:maxcliques2} (\ref{distance1}) that $|C \cap \Gamma(x)| = 0$ or $2$. In either case, $C$ contains vertices at distance at least $2$ from $x$. In particular, $C \cap \Gamma_2(x)$ is nonempty exactly when $d_\Gamma(x,C) = 1$ or $2$; in these cases the set $S$ defined by
	\begin{equation} \label{eq:S}
	S := \big\{ w \in \Gamma(x) \ : \ w \notin C \text{ and } w \in \mu(x,y) \text{ for some } y \in C \cap \Gamma_2(x) \big\}
	\end{equation}
has at least $4$ vertices.

\begin{lemma} \label{lemma:mu-clique}
Assume that $\Gamma$ is locally $n \times n$ grid with $n \geq 3$. Let $x \in \V(\Gamma)$ and $C$ a maximal clique in $\Gamma$ with $d_\Gamma(x,C) = 1$ or $2$. Let $\Delta$ be the union of all graphs $\mu(x,y)$ where $y \in C \cap \Gamma_2(x)$, and let $S$ be as in (\ref{eq:S}). Define
	\[ T = \left\{ \begin{aligned} &\varnothing &&\text{if } d_\Gamma(x,C) = 2, \\ &S \cap (\Gamma(u) \cup \Gamma(v)) &&\text{if } C \cap \Gamma(x) = \{u,v\}. \end{aligned} \right. \]
Then the following hold:
	\begin{enumerate}[(1)]
	\item \label{mu-number} Each vertex $w \in S$ lies in a unique $\mu$-graph in $\Delta$ if $w \in T$, and in exactly two $\mu$-graphs in $\Delta$ if $w \in S \setminus T$.
	\item \label{mu-meet} The intersection of any two distinct $\mu$-graphs in $\Delta$ is either $C \cap \Gamma(x)$, or the union of $C \cap \Gamma(x)$ and an edge in $S \setminus T$. 
	\item \label{matching} The set of all edges $e$ such that $e$ is in exactly two $\mu$-graphs in $\Delta$ is a perfect matching on $S \setminus T$.
	\end{enumerate}
\end{lemma}

\begin{proof}
We first prove statement (\ref{mu-number}). Let $w \in S$. Then $w \sim_\Gamma y$ for some $y \in C$, so that $d_\Gamma(w,C) = 1$. Thus $|\Gamma(w) \cap C| = 2$ by Lemma \ref{lemma:maxcliques2} (\ref{distance1}), and there is a unique vertex $z \in \Gamma(w) \cap C$ distinct from $y$. Hence $w$ lies in at most two $\mu$-graphs in $\Delta$, since each such $\mu$-graph is $\mu(x,y')$ for some $y' \in \Gamma(w) \cap C \cap \Gamma_2(x)$. If $w \notin T$ then $\Gamma(w) \cap C \cap \Gamma(x) = \varnothing$, so we must have $z \in \Gamma_2(x)$ and $w \in \mu(x,z)$. Otherwise $z \in \Gamma(x)$, so $\mu(x,y)$ is the unique $\mu$-graph in $\Delta$ that contains $w$. This proves statement (\ref{mu-number}).

We now prove statement (\ref{mu-meet}). First observe that $C \cap \Gamma(x)$ is contained in every $\mu$-graph in $\Delta$: this is vacuously true if $d_\Gamma(x,C) = 2$ since then $C \cap \Gamma(x) = \varnothing$, while if $d_\Gamma(x,C) = 1$ then $C \cap \Gamma(x) \subseteq \Gamma(x) \cap \Gamma(y) = \mu(x,y)$ for each $y \in C \setminus \Gamma(x)$, since $C$ is a clique. Let $\mu(x,y)$ and $\mu(x,z)$ be distinct $\mu$-graphs in $\Delta$ with a common vertex $w \notin C \cap \Gamma(x)$. Then $\Gamma(w) \cap C = \{y,z\} \subseteq \Gamma_2(x)$ by Lemma \ref{lemma:maxcliques2}, so $\Gamma(w) \cap C \cap \Gamma(x) = \varnothing$ and hence $w \in S \setminus T$. Also $\{w,y,z\}$ is a triangle, so by Lemma \ref{lemma:basic-nxn} (\ref{basic-triangles}) there is a unique maximal clique $C'$ of $\Gamma$ which contains $\{w,y,z\}$. Then $w \in C' \cap \Gamma(x)$, so $d_\Gamma(x,C') = 1$ and $|C' \cap \Gamma(x)| = 2$ by Lemma \ref{lemma:maxcliques2} (\ref{distance1}). Let $w'$ be the unique vertex in $C' \cap \Gamma(x)$ distinct from $w$. Then $w'$ is a common vertex of $\mu(x,y)$ and $\mu(x,z)$, and $w' \in \Gamma(w)$. So $w' \in S \setminus T$ by statement (\ref{mu-number}) above, as illustrated in Figure \ref{figure:mu-clique} on the left. Suppose that there is a third vertex $w''$ common to $\mu(x,y)$ and $\mu(x,z)$ with $w'' \notin C \cap \Gamma(x)$. Then again there is a unique maximal clique $C''$ of $\Gamma$ containing $\{w'',y,z\}$. Now $C' \neq C$, and it follows from Lemma \ref{lemma:basic-nxn} (\ref{basic-maxcliques}) that $C$ and $C'$ are the only two maximal cliques containing the edge $\{y,z\}$. So $C''$ is either $C$ or $C'$. Since $w'' \notin \{w,w'\} = C' \cap \Gamma(x)$, $C'' \neq C'$. Therefore $C'' = C$, so $w'' \in C \cap \Gamma(x)$, contradiction. Hence $\mu(x,y) \cap \mu(x,z) = (C \cap \Gamma(x)) \cup \{w,w'\}$, and $\{w,w'\}$ is an edge in $S \setminus T$. This proves statement (\ref{mu-meet}).

Finally we prove statement (\ref{matching}). Let $\Delta'$ denote the subgraph of $\Gamma$ consisting of all edges $e$ such that $e$ is in exactly two $\mu$-graphs in $\Delta$. By statement (\ref{mu-number}) above we have $\V(\Delta') \subseteq S \setminus T$. Also each $w \in S \setminus T$ is contained in exactly two $\mu$-graphs $\mu(x,y)$ and $\mu(x,z)$ in $\Delta$, and by statement (\ref{mu-meet}) we have $\mu(x,y) \cap \mu(x,z) = (C \cap \Gamma(x)) \cup e$ for some edge $e$ in $S \setminus T$. So $e = \{w,w'\}$ for some $w' \in S \setminus T$. It follows from statement (\ref{mu-number}) that $\mu(x,y)$ and $\mu(x,z)$ are the only $\mu$-graphs in $\Delta$ which contain $w'$, and thus $\{w,w'\} \in \E(\Delta')$ by the definition of $\Delta'$. So $w \in \V(\Delta')$ and $w$ is contained in an edge of $\Delta'$. Since $w$ is arbitrary, it follows that $S \setminus T \subseteq \V(\Delta')$, so $S \setminus T = \V(\Delta')$, and each vertex of $S \setminus T$ lies in an edge of $\Delta'$. Suppose that there are two distinct edges $e_1$ and $e_2$ of $\Delta'$ with a common vertex $w$. Then $e_1 = \mu(x,y) \cap \mu(x,z)$ and $e_2 = \mu(x,y') \cap \mu(x,z')$ for some $y, z, y', z' \in C \cap \Gamma_2(x)$ with $\{y,z\} \neq \{y',z'\}$. Thus $\{y,z,y',z'\} \subseteq \Gamma(w) \cap C$, as illustrated in Figure \ref{figure:mu-clique} on the right. This implies that $|\Gamma(w) \cap C| \geq 3$, contradiction. Therefore no two edges in $\Delta'$ have a common vertex, and so $\Delta'$ is a perfect matching on $S \setminus T$. This proves statement (\ref{matching}).
\end{proof}

\begin{center}
\begin{figure}
\begin{pspicture}(-6,-1.5)(6,2.5)
\rput(-4,0){
	\xdiagram
	\rput(-1.5,-0.5){\pnode(0,0.25){w} \pnode(0,-0.25){w'}} \rput(1.5,0){\pnode(0,0.25){y} \pnode(0,-0.25){z}}
	\ncline[linestyle=dashed]{w}{y} \ncline[linestyle=dashed]{w'}{z}
	\psellipse[linestyle=dashed,fillstyle=solid,fillcolor=white](-1.5,-0.5)(0.55,0.65)
	\psellipse[linestyle=dashed](1.5,0)(0.75,0.85)
	\rput{7}{\psellipse*[linecolor=white](-1.05,-0.25)(0.1,0.22)} \rput{9}{\psellipse*[linecolor=white](0.755,-0.27)(0.1,0.203)} \pscircle*[linecolor=white](0.9,-0.35){2pt} \psline[linewidth=1.1pt,linecolor=lightgray](0.7,0.09)(0.8,0.065) \psline[linewidth=1.1pt,linecolor=lightgray](0.7,0.115)(0.97,0.16)
	\psellipse*[linecolor=lightgray](1.5,0)(0.55,0.65)
	\rput(-1.5,-0.5){\cnode*(0,0.25){2pt}{w} \cnode*(0,-0.25){2pt}{w'} \ncline{w}{w'} \rput(-0.3,0.25){\small $w$} \rput(-0.25,-0.2){\small $w'$}}
	\rput(1.5,0){\cnode*(0,0.25){2pt}{y} \cnode*(0,-0.25){2pt}{z} \ncline{y}{z} \rput(0.25,0.25){\small $y$} \rput(0.25,-0.25){\small $z$}}
	\rput(0,-0.25){\small $C'$}
}
\rput(4,0){
	\xdiagram
	\psellipse*[linecolor=lightgray](1.5,0)(0.75,0.9)
	\rput(-1.2,-0.25){\cnode*(-0.9,0.25){2pt}{e1} \cnode*(0,0){2pt}{w} \cnode*(-0.25,-0.9){2pt}{e2} \ncline{e1}{w} \ncline{w}{e2}
	\rput(-0.45,-0.05){\small $e_1$} \rput(0.1,-0.55){\small $e_2$} \rput(0,0.25){\small $w$}}
	\rput(1.5,0){\cnode*(-0.2,0.5){2pt}{y} \cnode*(0.3,0.2){2pt}{z} \cnode*(0.3,-0.2){2pt}{y'} \cnode*(-0.2,-0.5){2pt}{z'}
	\rput(-0.2,0.8){\small $y$} \rput(0.3,0.45){\small $z$} \rput(0.3,-0.45){\small $y'$} \rput(-0.2,-0.8){\small $z'$}}
	\ncline{w}{y} \ncline{w}{z} \ncline{w}{y'} \ncline{w}{z'}
}
\end{pspicture}
\caption{$\Gamma$ and $C$ as in Lemma \ref{lemma:mu-clique}, with $C$ shown in gray and $\Gamma(x) \cap C = \varnothing$ if $d_\Gamma(x,C) = 2$}
\label{figure:mu-clique}
\end{figure}
\end{center}

\begin{corollary} \label{corollary:mu-clique}
Assume that $\Gamma$ is locally $n \times n$ grid with $n \geq 3$. Let $x$, $C$, and $S$ be as in Lemma \ref{lemma:mu-clique}. Then
	\[ 2|S| = \sum_{y \in C \cap \Gamma_2(x)} c_2(x,y) \equiv 0 \pmod{4}. \]
\end{corollary}

\begin{proof}
We count in two ways the number $\sigma := \big|\{ (w,y) \ : \ y \in C \cap \Gamma_2(x), \ w \in \Gamma(x) \cap \Gamma(y) \}\big|$. First we have
	\[ \sigma = \sum_{y \in C \cap \Gamma_2(x)} \big|\{ w \in \Gamma(x) \ : \ w \sim_\Gamma y \}\big| = \sum_{y \in C \cap \Gamma_2(x)} |\mu(x,y)|. \]
Next we have
	\begin{equation}
	\sigma
	= \sum_{w \in \Gamma(x)} \big|\{ y \in C \cap \Gamma_2(x) \ : \ y \sim_\Gamma w \}\big|
	= \sum_{w \in \Gamma(x)} \big|\{ y \in C \cap \Gamma_2(x) \ : \ w \in \mu(x,y) \}\big|. \label{eq:sigma}
	\end{equation}
The nonzero terms in the sum on the right side of (\ref{eq:sigma}) correspond exactly to those $w \in \V(\Delta) = S \cup (C \cap \Gamma(x))$. We apply Lemma \ref{lemma:mu-clique} to the right side of (\ref{eq:sigma}). If $d_\Gamma(x,C) = 2$ then $C \cap \Gamma(x) = \varnothing$, and each vertex in $S$ lies in exactly two $\mu$-graphs in $\Delta$. Hence
	\[ \sum_{w \in \Gamma(x)} \big|\{ y \in C \cap \Gamma_2(x) \ : \ w \in \mu(x,y) \}\big| = \sum_{w \in S} \big|\{ y \in C \cap \Gamma_2(x) \ : \ w \in \mu(x,y) \}\big| = 2|S|. \]
Now suppose that $d_\Gamma(x,C) = 1$. Then $|C \cap \Gamma(x)| = 2$; let $\{u,v\} = C \cap \Gamma(x)$. Let $T$ denote the set of all vertices $w \in S$ such that $\Gamma(w) \cap C \cap \Gamma(x) \neq \varnothing$, that is, $T = S \cap (\Gamma(u) \cup \Gamma(v))$. By Lemma \ref{lemma:mu-clique} (\ref{mu-number}) each vertex in $S \setminus T$ is in exactly two $\mu$-graphs in $\Delta$, while each vertex in $T$ is in a unique $\mu$-graph in $\Delta$. Note also that no vertex in $T$ is adjacent to both $u$ and $v$, for otherwise such a vertex will have three neighbours in $C$, contradiction. It follows that $S \cap \Gamma(u) \cap \Gamma(v) = \varnothing$ and $T$ is the disjoint union of $S \cap \Gamma(u)$ and $S \cap \Gamma(v)$, as illustrated in Figure \ref{figure:mu-vertices-cor}. By Lemma \ref{lemma:maxcliques2} (\ref{distance1}) we have $|C \cap \Gamma_2(x)| = n-1$, so there are $n-1$ $\mu$-graphs in $\Delta$, each of which contains $\{u,v\}$. Since each vertex in $T$ lies in a unique $\mu$-graph in $\Delta$, it follows that $|S \cap \Gamma(u)| = |S \cap \Gamma(v)| = n-1$. Thus $|T| = 2(n-1)$. Recalling that in the right side of (\ref{eq:sigma}), the nonzero contributions come from $w \in S \cup \{u,v\}$, we have
	\begin{align*}
	\sum_{w \in \Gamma(x)} \big|\{ y \in C \cap \Gamma_2(x) \ : \ y \sim_\Gamma w \}\big|
	&= \sum_{w \in S \cup \{u,v\}} \big|\{ y \in C \cap \Gamma_2(x) \ : \ y \sim_\Gamma w \}\big| \\
	&= 2|S \setminus T| + (n-1) |\{u,v\}| + 1 \cdot |T| \\
	&= 2(|S| - 2(n-1)) + (n-1) \cdot 2 + 1 \cdot 2(n-1) \\
	&= 2|S|.
	\end{align*}
Therefore in both cases $2|S| = \sigma$. By Lemma \ref{lemma:mu-clique} (\ref{matching}) the set $S \setminus T$ can be paritioned into subsets of size $2$, so $|S \setminus T|$ must be even. Since also $|T| = 2(n-1)$ is even, it follows that $|S|$ is even. Hence $\sigma = 2|S| \equiv 0 \pmod{4}$, which completes the proof.
\end{proof}

\begin{center}
\begin{figure}
\begin{pspicture}(-2.7,-1.5)(2.7,2.5)
\psellipse(0,0)(1,1.5) \rput(0,1.8){\small $\Gamma(x)$}
\cnode*(-0.2,0.9){2pt}{u} \cnode*(0.2,0.9){2pt}{v} \rput(-0.2,1.15){\small $u$} \rput(0.2,1.15){\small $v$} \ncline{u}{v}
\pnode(-0.5,0.55){Tu1} \pnode(0.5,0.55){Tv1} \ncline{u}{Tu1} \ncline{v}{Tv1}
\pnode(-0.5,-0.25){Tu2} \pnode(0.5,-0.25){Tv2} \pnode(-0.35,-0.9){Su} \pnode(0.35,-0.9){Sv} \ncline{Tu2}{Su} \ncline{Tv2}{Sv}
\psellipse[fillcolor=white,fillstyle=solid](-0.5,0.15)(0.25,0.4) \psline(-0.65,0.15)(-1.25,0.15) \rput(-2,0.15){\small $S \cap \Gamma(u)$}
\psellipse[fillcolor=white,fillstyle=solid](0.5,0.15)(0.25,0.4) \psline(0.65,0.15)(1.25,0.15) \rput(2,0.15){\small $S \cap \Gamma(v)$}
\rput(0,-0.9){\psellipse[fillcolor=white,fillstyle=solid](0,0)(0.55,0.35) \rput(0,0){\small $S \setminus T$}}
\end{pspicture}
\caption{$\V(\Delta) = S \cup \{u,v\}$, $T = (S \cap \Gamma(u)) \cup (S \cap \Gamma(v))$} \label{figure:mu-vertices-cor}
\end{figure}
\end{center}

For the next result we assume that all $\mu$-graphs have order at least $2(n-1)$. In this case Lemma \ref{lemma:maxcliques2} states that any $(x,y) \in \Gamma_2$ satisfies the following: if $c_2(x,y) = 2n$ then $d_\Gamma(x,C) = 1$ for all maximal cliques $C$ of $\Gamma$ containing $y$, and if $c_2(x,y) = 2(n-1)$ then $d_\Gamma(x,C) = 1$ for $2(n-1)$ maximal cliques $C$ containing $y$ and $d_\Gamma(x,C) = 2$ for the remaining two cliques.

\begin{lemma} \label{lemma:cliques-bounded}
Assume that $\Gamma$ is connected and locally $n \times n$ grid, and that all $\mu$-graphs in $\Gamma$ have order at least $2(n-1)$. Then $\diam(\Gamma) \leq 3$, and $d_\Gamma(x,C) = 1$ or $2$ for any $x \in \V(\Gamma)$ and maximal clique $C$ not containing $x$. Furthermore, the following hold:
	\begin{enumerate}[(1)]
	\item \label{distance1-bounded} If $d_\Gamma(x,C) = 1$ then the number of vertices $y \in C \cap \Gamma_2(x)$ satisfying $c_2(x,y) = 2(n-1)$ is even. Moreover, if $n$ is even then $c_2(x,z) = 2n$ for some $z \in C \cap \Gamma_2(x)$.
	\item \label{distance2-bounded} If $d_\Gamma(x,C) = 2$ then each vertex $y \in C \cap \Gamma_2(x)$ satisfies $c_2(x,y) = 2(n-1)$. Moreover, if $n$ is even then $C \nsubseteq \Gamma_2(x)$.
	\end{enumerate}
\end{lemma}

\begin{proof}
Let $x \in \V(\Gamma)$ and let $C$ be a maximal clique in $\Gamma$ not containing $x$. It follows from Lemma \ref{lemma:parameters} (\ref{b3,c3}) that $b_3(x,y) \leq (n-(n-1)-1)^2 = 0$ for all $y \in \Gamma_3(x)$, which implies that $k_i(x) = 0$ for all $i \geq 4$. By connectedness $\diam(\Gamma) \leq 3$. Using $m = n-1$ in (\ref{eq:k3}) we have
	\[ |\Gamma_3(x)| = k_3(x) \leq \frac{n-1}{2} < n+1 = |C|, \]
and thus $C \nsubseteq \Gamma_3(x)$. Since $\diam(\Gamma) \leq 3$, either $d_\Gamma(x,C) = 1$ or $d_\Gamma(x,C) = 2$. Let $r = |C \cap \Gamma_2(x)|$ and let $s$ be the number of vertices $y \in C \cap \Gamma_2(x)$ with $c_2(x,y) = 2(n-1)$. Then $c_2(x,y) = 2n$ for the remaining $r-s$ vertices, and
	\[ \sum_{y \in C \cap \Gamma_2(x)} c_2(x,y) = s \cdot 2(n-1) + (r-s) \cdot 2n = 2(rn - s). \]
By Corollary \ref{corollary:mu-clique} this number is divisible by $4$. Hence $t := rn - s$ is even. Suppose first that $d_\Gamma(x,C) = 1$. Then $r = n-1$ by Lemma \ref{lemma:maxcliques2} (\ref{distance1}), so $t = n(n-1) - s$, and since $t$ is even, $s$ must also be even. This proves the first part of statement (\ref{distance1-bounded}). If in this case $n$ is even then $r = n-1$ is odd and so $r \neq s$ since $s$ is even, whence $c_2(x,z) = 2n$ for some $z \in C \cap \Gamma_2(x)$ and statement (\ref{distance1-bounded}) is proved. Now suppose that $d_\Gamma(x,C) = 2$. Then $c_2(x,y) \leq 2(n-1)$ for any $y \in C \cap \Gamma_2(x)$ by Lemma \ref{lemma:maxcliques2} (\ref{distance2}), and together with the hypothesis this gives us $c_2(x,y) = 2(n-1)$ for all $y \in C \cap \Gamma_2(x)$. So the first part of statement (\ref{distance2-bounded}) holds. Thus $s = r$; also $r \geq (n-1) + 1 = n$ by the second part of Lemma \ref{lemma:maxcliques2} (\ref{distance2}), so $r = n$ or $n+1$. If $r = n$ then $t = n^2 - n$, which is even for any $n$. If $r = n+1$ then $t = n^2 - 1$, which is even exactly when $n$ is odd; it follows that if $n$ is even then $|C \cap \Gamma_2(x)| = r \neq n+1 = |C|$, and thus $C \nsubseteq \Gamma_2(x)$. This completes the proof of statement (\ref{distance2-bounded}).
\end{proof}

%% file: locallygrid-sec6-lowerbound.tex
\section{Graphs with $|\mu| \geq 2(n-1)$} \label{sec:lowerbound}

In this section we consider the special case where all $\mu$-graphs have order at least $2(n-1)$. The main result is Theorem \ref{maintheorem:order-diameter}. We apply these results to locally $3 \times 3$ grid and locally $5 \times 5$ grid graphs.

The subcase where all $\mu$-graphs have the maximum possible order $2n$ is covered by the remarks following \cite[Lemma, Section 2]{4by4}. We state below this result for locally $n \times n$ grid graphs.

\begin{theorem} \label{theorem:srg} \cite[Section 2]{4by4}
Assume that $\Gamma$ is connected and locally $n \times n$ grid, and that all $\mu$-graphs of $\Gamma$ have order $2n$. Then $\diam(\Gamma) = 2$ and $\Gamma$ is strongly regular with parameters
	\[ \left( \frac{n^3 + n + 2}{2}, \ n^2, \ 2(n-1), \ 2n \right). \]
\end{theorem}

Indeed, in this subcase equation (\ref{eq:b=c}) with $i = 2$ and $i = 3$ gives us
	\[ k_2(x) = \frac{n(n-1)^2}{2} \ \ \text{and} \ \ k_3(x) = 0 \]
for all $x \in \V(\Gamma)$, $\diam(\Gamma) = 2$. All $\mu$-graphs have the same order so $c_2(x,y)$ is constant for all $\{x,y\} \in \Gamma_2$; this together with Lemma \ref{lemma:basic-nxn} (\ref{basic-edges}) imply that $\Gamma$ is strongly regular.

Suppose now that some $\mu$-graph in $\Gamma$ has order $2(n-1)$. By (\ref{eq:k2}) and (\ref{eq:k3})
	\begin{equation} \label{eq:k2,k3}
	k_2(x) \leq \frac{n^2(n-1)}{2} \ \ \text{and} \ \ k_3(x) \leq \frac{n-1}{2}
	\end{equation}
for all $x \in \V(\Gamma)$. With $k_{2,2(n-1)}(x)$ and $k_{2,2n}(x)$ as in (\ref{eq:k_2,2m}), counting the number of edges between $\Gamma(x)$ and $\Gamma_2(x)$ yields the following special case of (\ref{eq:sum-k_2,2m}):
	\begin{equation} \label{eq:edges-k2}
	n^2(n-1)^2 = 2(n-1)\,k_{2,2(n-1)}(x) + 2n\,k_{2,2n}(x).
	\end{equation}
The left side of (\ref{eq:edges-k2}) is divisible by $2n(n-1)$. Hence $k_{2,2(n-1)}(x) \equiv 0 \pmod{n}$, $k_{2,2n}(x) \equiv 0 \pmod{n-1}$, and
	\[ \frac{n(n-1)}{2} = \frac{k_{2,2(n-1)}(x)}{n} + \frac{k_{2,2n}(x)}{n-1}. \]
By (\ref{eq:k2-sum}) we have $k_2(x) = k_{2,2(n-1)}(x) + k_{2,2n}(x)$, and substituting from this into the equation above gives us
	\begin{equation} \label{eq:k2-sum2}
	k_2(x) = \frac{n(n-1)^2}{2} + \frac{k_{2,2(n-1)}(x)}{n}.
	\end{equation}

\begin{proof}[Proof of Theorem \ref{maintheorem:order-diameter}]
Assume that all $\mu$-graphs have order at least $2(n-1)$. By Lemma \ref{lemma:cliques-bounded}, $\diam(\Gamma) \leq 3$. Thus $|\V(\Gamma)| = 1 + n^2 + k_2(x) + k_3(x)$ for any $x \in \V(\Gamma)$; this together with (\ref{eq:k2,k3}) gives the bound
	\[ |\V(\Gamma)| \leq 1 + n^2 + \frac{n^2(n-1)}{2} + \frac{n-1}{2} = \frac{(n^2 + 1)(n+1)}{2} \]
and hence
	\begin{equation} \label{eq:order-bounded}
	|\V(\Gamma)| \leq \left\lfloor\frac{(n^2 + 1)(n+1)}{2}\right\rfloor.
	\end{equation}
This proves Theorem \ref{maintheorem:order-diameter} (\ref{>2(n-1)}).

We claim that for any $x \in \V(\Gamma)$, $k_2(x) = n^2(n-1)/2$ if and only if $c_2(x,y) = 2(n-1)$ for all $y \in \Gamma_2(x)$. Indeed, if $c_2(x,y) = 2(n-1)$ for all $y \in \Gamma_2(x)$ then $k_{2,2(n-1)}(x) = k_2(x)$ and $k_{2,2m}(x) = 0$ for all $m \neq n-1$, and it follows from (\ref{eq:k2-sum2}) that
	\[ k_2(x) = \frac{n(n-1)^2}{2} \cdot \frac{n}{n-1} = \frac{n^2(n-1)}{2}. \]
Conversely, suppose that $k_2(x) = n^2(n-1)/2$. Then from (\ref{eq:k2-sum2}) we get
	\begin{align*}
	k_{2,2(n-1)}(x)
	&= n\,\left(k_2(x) - \frac{n(n-1)^2}{2}\right) \\
	&= n\,\left(\frac{n^2(n-1)}{2} - \frac{n(n-1)^2}{2}\right) \\
	&= \frac{n^2(n-1)}{2} \\
	&= k_2(x).
	\end{align*}
So $c_2(x,y) = 2(n-1)$ for all $y \in \Gamma_2(x)$, which proves the claim.

We now prove Theorem \ref{maintheorem:order-diameter} (\ref{=2(n-1)}). Assume first that equality holds in (\ref{eq:order-bounded}). Let $x \in \V(\Gamma)$. It follows from (\ref{eq:k2,k3}) that $k_2(x) = n^2(n-1)/2$. Thus, by the claim, $c_2(x,y) = 2(n-1)$ for all $y \in \Gamma_2(x)$. Since $x$ is arbitrary this proves that all $\mu$-graphs have order $2(n-1)$.

Conversely, assume that all $\mu$-graphs have order $2(n-1)$. Let $x \in \V(\Gamma)$ be arbitrary. Then $c_2(x,y) = 2(n-1)$ for all $y \in \Gamma_2(x)$, and thus $k_2(x) = n^2(n-1)/2$ by the claim. Suppose that $k_3(x) = 0$. Then
	\[ |\V(\Gamma) = 1 + n^2 + \frac{n^2(n-1)}{2} = \frac{n^3 + n^2 + 2}{2}. \]
By Lemma \ref{lemma:basic-nxn} (\ref{basic-triangles}), $n+1$ divides $2|\V(\Gamma)| = n^3 + n^2 + 2$, and hence $n+1$ divides $2$, contradiction. Thus $k_3(x) \neq 0$ and $\diam(\Gamma) = 3$. By (\ref{eq:k2,k3}) we have $2\,k_3(x) \leq n-1$, so that $2\,k_3(x) = n-s$ for some $s \in \{1, \ldots, n-1\}$. Therefore
	\[ 2|\V(\Gamma)| = 2 + 2n^2 + 2\,k_2(x) + 2\,k_3(x) = n^3 + n^2 + n - s + 2, \]
which is divisible by $n+1$ if and only if $s = 1$. So $k_3(x) = (n-1)/2$ and $n$ is odd, and consequently equality holds in (\ref{eq:order-bounded}). This proves the first part of Theorem \ref{maintheorem:order-diameter}, and also that in case of equality we have $\diam(\Gamma) = 3$ and $n$ odd.

It remains to show that $\Gamma$ is a distance-regular antipodal cover of $K_{n^2+1}$ whenever all $\mu$-graphs have order $2(n-1)$. By the hypothesis $c_2$ is constant on $\Gamma_2$, so we need only to show that $b_2$ and $c_3$ are constant on $\Gamma_2$ and $\Gamma_3$, respectively. Let $x \in \V(\Gamma)$. By Lemma \ref{lemma:parameters} (\ref{b3,c3}) any $y \in \Gamma_3(x)$ satisfies $c_3(x,y) \geq ((n-1)+1)^2 = n^2 = |\Gamma(y)| \geq c_3(x,y)$, so $c_3(x,y) = n^2$. Thus $c_3$ is constant for any pair of vertices in $\Gamma_3$. By Lemma \ref{lemma:parameters} (\ref{b2}) any $y \in \Gamma_2(x)$ satisfies $b_2(x,y) \leq (n-(n-1))^2 = 1$. Letting $b_{2,1}(x) = \big|\{ y \in \Gamma_2(x) \ : \ b_2(x,y) = 1 \}\big|$ and applying (\ref{eq:b=c}) with $i = 3$ we get
	\[ b_{2,1}(x) = \sum_{y \in \Gamma_2(x)} b_2(x,y) = \sum_{z \in \Gamma_3(x)} c_3(x,y) = n^2 k_3(x) = \frac{n^2(n-1)}{2} = k_2(x). \]
Hence $b_2(x,y)= 1$ for all $y \in \Gamma_2(x)$, which shows that $b_2$ is constant on $\Gamma_2$. Therefore $\Gamma$ is distance-regular. For all $z \in \Gamma_3(x)$ we have $a_3(x,z) = n^2 - c_3(x,z) = 0$, so no two vertices $w, z \in \Gamma_3(x)$ are adjacent in $\Gamma$. If $d_\Gamma(w,z) = 2$ then there is a vertex $y \in \Gamma_2(x)$ such that $w, z \in \Gamma(y)$; in this case $b_2(x,y) > 1$, contradiction. So $d_\Gamma(w,z) = 3$ for any $w, z \in \Gamma_3(x)$. Therefore $x \cup \Gamma_3(x)$ is an antipodal block for all $x \in \V(\Gamma)$, and the quotient graph with respect to the resulting partition is $K_{n^2+1}$. This completes the proof of Theorem \ref{maintheorem:order-diameter} (\ref{=2(n-1)}).
\end{proof}

In the remainder of this section we apply the above results to locally $n \times n$ grid graphs for $n \in \{3,5\}$. We will use the following technical lemma:

\begin{lemma} \label{lemma:ell}
Assume that $\Gamma$ is locally $n \times n$ grid, and that all $\mu$-graphs in $\Gamma$ have order at least $2(n-1)$. Let $x \in \V(\Gamma)$ and $k_{2,2(n-1)}(x)$ as in (\ref{eq:k_2,2m}). Then $k_{2,2(n-1)}(x) = \ell_x n$, for some integer $\ell_x \leq n(n-1)/2$ such that
	\[ \ell_x + k_3(x) \equiv \left\{\begin{aligned} &0 \pmod{n+1} &&\text{if $n$ is even}; \\ &0 \pmod{(n+1)/2} &&\text{if $n$ is odd}. \end{aligned}\right. \]
\end{lemma}

\begin{proof}
Recall from the remarks after Theorem \ref{theorem:srg} that $k_{2,2(n-1)}(x) \equiv 0 \pmod{n}$, so indeed $k_{2,2(n-1)}(x) = \ell_x n$ for some $\ell_x$. From (\ref{eq:k2,k3}) and the definition of $k_{2,2(n-1)}(x)$ we have
	\[ k_{2,2(n-1)}(x) \leq k_2(x) \leq \frac{n^2(n-1)}{2}, \]
and hence $\ell_x \leq n(n-1)/2$. By (\ref{eq:k2-sum2}) we have $k_2(x) = n(n-1)^2/2 + \ell_x$, and since $\diam(\Gamma) \leq 3$ by Theorem \ref{maintheorem:order-diameter} (\ref{>2(n-1)}),
	\[ |\V(\Gamma)| = 1 + n^2 + k_2(x) + k_3(x) = \frac{(n+1)(n^2 - n + 2)}{2} + \ell_x + k_3(x). \]
Recall from Lemma \ref{lemma:basic-nxn} (\ref{basic-triangles}) that $n+1$ divides $2|\V(\Gamma)|$. So $2(\ell_x + k_3(x)) \equiv 0 \pmod{n+1}$, and the result follows.
\end{proof}

\subsection{The subcase $n = 3$}

If $\Gamma$ is locally $3 \times 3$ grid then by Lemma \ref{lemma:basic-nxn} (\ref{basic-mu}) any $\mu$-graph of $\Gamma$ has order at least $4 = 2(n-1)$. Hence Theorem \ref{maintheorem:order-diameter} (\ref{>2(n-1)}) may be applied. The locally $3 \times 3$ grid graphs belong to a more general family classified by Hall in \cite{3byq}.

\begin{proposition} \label{proposition:3x3}
Assume that $\Gamma$ is connected and locally $3 \times 3$ grid. Then all $\mu$-graphs of $\Gamma$ have the same order $|\mu| \in \{4,6\}$. Moreover:
	\begin{enumerate}[1.]
	\item \label{mu=4} If $|\mu| = 4$ then $\Gamma \cong J(6,3)$ (equivalently, to the graph in Construction \ref{example:drg} with $n = 3$).
	\item \label{mu=6} If $|\mu| = 6$ then $\Gamma$ is isomorphic to the complement $\overline{K_4 \square K_4}$ of the $4 \times 4$ grid graph.
	\end{enumerate}
\end{proposition}

\begin{remark} \label{remark:3x3}
Proposition \ref{proposition:3x3} follows from more general results \cite[Theorems 1 and 2]{3byq} concerning line graphs of certain partial linear spaces of order $2$, which are locally $3 \times n$ grid for some $n$. We give a self-contained elementary proof of the subclass of locally $3 \times 3$ grid graphs based on the theory developed in our paper. We note that the two examples we obtain in Proposition \ref{proposition:3x3} all come from partial linear spaces $\mathscr{T}(\Omega,\Omega')$ in \cite[Theorem 1]{3byq}, in particular, $J(6,3)$ arises from $|\Omega| = 6$, $\Omega' = \varnothing$; and $\overline{K_4 \square K_4}$ arises in two ways, namely, $(|\Omega|,|\Omega'|) = (4,1)$ or $(3,2)$. The graph $J(6,3)$ also arises, for example, from the space $\mathscr{S}p(V,f)$ where $f$ is nondegenerate and $V = \F_2^4$.
\end{remark}

\begin{proof}[Proof of Proposition \ref{proposition:3x3}]
Assume that $\Gamma$ is locally $3 \times 3$ grid. By Theorem \ref{maintheorem:order-diameter} (\ref{>2(n-1)}), we have $\diam(\Gamma) \leq 3$ and $|\V(\Gamma)| \leq 20$. Applying (\ref{eq:k2,k3}), we obtain for any vertex $x$ that $k_2(x) \leq 9$ and $k_3(x) \leq 1$. Also, by Lemma \ref{lemma:parameters} (\ref{b2}) and \ref{lemma:parameters} (\ref{b3,c3}), $b_2(x,y) \leq 1$ for any $y \in \Gamma_2(x)$, and $c_3(x,y) \geq 9$ for any $y \in \Gamma_3(x)$.

We claim that if $\epsilon(x) = 2$, where $\epsilon$ is as in (\ref{eq:epsilon}), then $\big(k_{2,4}(x), \, k_{2,6}(x)\big)$ is either $(0,6)$ or $(6,2)$. Indeed, by Lemma \ref{lemma:ell} we have $k_{2,4}(x) = 3\ell_x$ for some integer $\ell_x$ satisfying $\ell_x \leq 3$ and $\ell_x + k_3(x) \equiv 0 \pmod{2}$. Now $k_3(x) = 0$, so $\ell_x \in \{0,2\}$, and from (\ref{eq:edges-k2}),
	\[ k_{2,6}(x) = 6 - \frac{2}{3}\, k_{2,4}(x) = 6 - 2\ell_x. \]
The claim follows.

We consider two cases.

\emph{Case 1.} Suppose that $\diam(\Gamma) = 2$. Then $k_2(x) = |\V(\Gamma)| - 10$ for any $x \in \V(\Gamma)$, so that $k_2(x)$ is constant. Also $k_3(x) = 0$, so $\epsilon(x) = 2$, and thus by the claim above $\big(k_{2,4}(x), \, k_{2,6}(x)\big) \in \{(0,6), (6,2)\}$.

Suppose first that $k_{2,4}(x) = 6$ and $k_{2,6}(x) = 2$. Then we can denote the elements of $\Gamma_2(x)$ by $y_i$ ($1 \leq i \leq 6$) and $z_j$ ($1 \leq j \leq 2$), where $c_2(x,y_i) = 4$ and $c_2(x,z_j) = 6$ for each $i$ and each $j$. Let $S_i = \Gamma_2(x) \cap \Gamma(y_i)$ and $T_j = \Gamma_2(x) \cap \Gamma(z_j)$. Notice that $[S_i]$ and $[T_j]$ are subgraphs of $[\Gamma(y_i)]$ and $[\Gamma(z_j)]$, respectively, where $[\Gamma(y_i)] \cong [\Gamma(z_j)] \cong K_3 \square K_3$. By Lemma \ref{lemma:basic-nxn} (\ref{basic-mu}), any $4$-cycle in $K_3 \square K_3$ has two edges in two distinct vertical cliques and two edges in two distinct horizontal cliques, and so its complement consists of one vertical and one horizontal clique in $K_3 \square K_3$. So for each $i \in \{1, \ldots, 5\}$, each induced subgraph $[S_i]$ has five vertices and is isomorphic to $\overline{C_4 \cup K_1}$, as illustrated in Figure \ref{figure:3by3}. Likewise, any $6$-cycle in $K_3 \square K_3$ has three edges in three distinct horizontal cliques and three edges in three distinct vertical cliques, so each clique in $K_3 \square K_3$ contains two vertices of the $6$-cycle. It follows that its complement consists of three vertices no two of which belong in the same clique, that is, no two of which are adjacent. Thus each $[T_j]$ is an empty graph of order three, $3K_1$. Each $[S_i \cup \{y_i\}]$ has two vertices of valency $5$ (including $y_i$) and four vertices of valency $3$, and the neighbourhood in $[S_i \cup \{y_i\}]$ of any of these vertices contains an edge. It follows that for all $i$ and all $j$ we have $z_j \notin S_i$, which implies that $z_j \not\sim_\Gamma y_i$. Hence $T_1, T_2 \subseteq \{z_1,\,z_2\}$, a contradiction since $T_1$ and $T_2$ have three elements each.

It follows that $k_{2,4}(x) = 0$ and $k_{2,6}(x) = 6$. Hence $c_2(x,y) = 6$ for all $y \in \Gamma_2(x)$, and since $x$ is arbitrary this holds for all pairs $(x,y) \in \Gamma_2$. Therefore all $\mu$-graphs of $\Gamma$ have size $|\mu| = 6$. By Theorem \ref{theorem:srg}, $\Gamma$ is strongly regular with parameters $(16, 9, 4, 6)$. Up to isomorphism there are exactly two such graphs (\cite{srgdatabase} or \cite[Table 1.1, line 3]{srgBM}), and of these only $\overline{K_4 \square K_4}$ is locally $3 \times 3$ grid. Thus part (\ref{mu=6}) of the statement holds.

\emph{Case 2.} Suppose now that $\diam(\Gamma) = 3$. We show that $k_3(x) \neq 0$ for all $x \in \V(\Gamma)$. Indeed, $\diam(\Gamma) = 3$ implies that there exists $x \in \V(\Gamma)$ such that $k_3(x) \neq 0$. Then $k_3(x) = 1$, so that $\Gamma(z) \subseteq \Gamma_2(x)$ for the unique $z \in \Gamma_3(x)$. Thus $k_2(x) \geq |\Gamma(z)| = 9$. If $k_3(y) = 0$ for some vertex $y$ then $\epsilon(y) = 2$, and it follows from the claim above that $k_2(y) = k_{2,4}(y) + k_{2,6}(y) \in \{6, 8\}$. Hence $|\V(\Gamma)| \in \{16, 18\}$ and $9 \leq k_2(x) = |\V(\Gamma)| - (1 + k_1(x) + k_3(x)) = |\V(\Gamma)| - 11 \in \{5, 7\}$, contradiction. Therefore no such $y$ exists, and $k_3(x) \neq 0$ for any $x \in \V(\Gamma)$.

It follows that any $x \in \V(\Gamma)$ satisfies $k_3(x) = 1$, say $\Gamma_3(x) = \{z\}$, so that $\Gamma_2(x) \supseteq \Gamma(z)$ and $k_2(x) \geq 9$. Also, by equation (\ref{eq:k2-sum2}), we have $k_2(x) = 6 + \ell_x$, and as $\ell_x \leq 3$ we conclude that $\ell_x = 3$ and $k_2(x) = 9$. Hence $k_{2,4}(x) = 3\ell_x = 9$ and $k_{2,6}(x) = 0$, so $c_2(x,y) = 4$ for all $y \in \Gamma_2(x)$. Since $x$ is arbitrary this holds for all pairs $(x,y) \in \Gamma_2$. Therefore all $\mu$-graphs of $\Gamma$ have size $|\mu| = 4$. By Theorem \ref{maintheorem:order-diameter} (\ref{=2(n-1)}), $\Gamma$ is a distance-regular antipodal double cover of $K_{10}$, and hence has $20$ vertices. Applying \cite[Theorem 1]{4by4} we conclude that $\Gamma$ is the Johnson graph $J(6,3)$ as in part (\ref{mu=4}) of the statement.
\end{proof}

\begin{center}
\begin{figure}
\begin{pspicture}(-7,-2)(7,2)
\rput(-4,0){
\gridgraph3
\pspolygon[linewidth=1pt](-1,1)(0,1)(0,0)(-1,0)(-1,1)
\psarc[linestyle=dashed,linewidth=1pt](0,-4.732){3.864}{75}{105} \psline[linestyle=dashed,linewidth=1pt](-1,-1)(1,-1)
\psarc[linestyle=dashed,linewidth=1pt](4.732,0){3.864}{165}{195} \psline[linestyle=dashed,linewidth=1pt](1,-1)(1,1)
\multido{\n=-1+1}{3}{\rput(\n,0){\multido{\n=-1+1}{3}{\qdisk(0,\n){2pt}}}}
\rput(-0.5,1.75){\small $\mu(x,y_i)$} \psline(-0.5,1.5)(-0.5,1)
\rput(0.5,-1.5){\small $[S_i] \cong \overline{C_4 \cup K_1}$} \psline(0.5,-1)(0.5,-1.25)
\rput(-2,0){\small $\Gamma(y_i)$}
}

\multido{\n=0+60}{6}{\rput{\n}(0,0){\qdisk(0.87,0.5){2pt} \psline(0.87,-0.5)(0.87,0.5)}}
\psline(0,1)(-0.87,-0.5) \psline(0,1)(0,-1) \psline(0,1)(0.87,-0.5) \psline(-0.87,0.5)(0,-1) \psline(0.87,0.5)(0,-1)
\rput(0,1.25){\small $y_i$} \rput(0,-1.5){\small $[S_i \cup \{y_i\}]$}

\rput(4,0){
\gridgraph3
\psline[linewidth=1pt](-1,1)(0,1)(0,0)(1,0)(1,-1)
\psarc[linewidth=1pt](0,-4.732){3.864}{75}{105}
\psarc[linewidth=1pt](2.732,0){3.864}{165}{195}
\multido{\n=-1+1}{3}{\rput(\n,0){\multido{\n=-1+1}{3}{\qdisk(0,\n){2pt}}}}
\rput(2,0){\small $\Gamma(z_j)$}
\rput(-0.5,1.75){\small $\mu(x,z_j)$} \psline(-0.5,1.5)(-0.5,1)
\rput(0,-1.5){\small $[T_j] \cong \overline{K_3}$}
}
\end{pspicture}
\caption{$S_i$ and $T_j$ as in Case 1 of the proof of Proposition \ref{proposition:3x3}} \label{figure:3by3}
\end{figure}
\end{center}

\subsection{The subcase $n = 5$}

If $\Gamma$ is connected and locally $5 \times 5$ grid then $\diam(\Gamma) \leq 5$ by Theorem \ref{theorem:diameter}. In Lemma \ref{lemma:order-bound} we use this to obtain an upper bound for the order of $\Gamma$.

\begin{lemma} \label{lemma:order-bound}
Assume that $\Gamma$ is connected and locally $5 \times 5$ grid. Then either $\diam(\Gamma) = 5$ and $\Gamma \cong J(10,5)$, or $\diam(\Gamma) = 4$, $|\V(\Gamma)| \leq 270$, and $|\V(\Gamma)| \equiv 0 \pmod{6}$.
\end{lemma}

\begin{proof}
The first part follows immediately from Theorem \ref{theorem:diameter}. Suppose that $\diam(\Gamma) \leq 4$. Applying (\ref{eq:k2}), (\ref{eq:k3}), and (\ref{eq:ki}) yields
	\[ \quad k_2 \leq \frac{5^2 \cdot 4^2}{4} = 100, \quad k_3 \quad \leq \frac{100 \cdot 9}{9} = 100, \quad k_4 \quad \leq \left\lfloor \frac{100 \cdot 4}{9} \right\rfloor = 44. \]
Hence $|\V(\Gamma)| \leq 1 + 25 + k_2(x) + k_3(x) + k_4(x) = 270$. Now $5 \equiv 2 \pmod{3}$, so $3$ divides $|\V(\Gamma)|$ by Lemma \ref{lemma:basic-nxn} (\ref{basic-triangles}), and $|\V(\Gamma)|$ is even since $\Gamma$ has odd valency. Thus $|\V(\Gamma)| \equiv 0 \pmod{6}$, and the lemma is proved.
\end{proof}

\begin{lemma} \label{lemma:mu-constant}
Assume that $\Gamma$ is connected and locally $5 \times 5$ grid.
	\begin{enumerate}[(1)]
	\item \label{mu=8} If all $\mu$-graphs in $\Gamma$ have constant order $|\mu|$, then either $|\mu| = 4$ and $\Gamma \cong J(10,5)$ or $\frac{1}{2}J(10,5)$, or $|\mu| = 8$.
	\item \label{mu>8-lem} If all $\mu$-graphs in $\Gamma$ have order at least $8$, then $\diam(\Gamma) \leq 3$ and $|\V(\Gamma)| \leq 78$.
	\end{enumerate}
\end{lemma}

\begin{proof}
Statement \eqref{mu>8-lem} follows immediately from Theorem \ref{maintheorem:order-diameter}.

Assume that all $\mu$-graphs in $\Gamma$ have constant order $|\mu|$. By Lemma \ref{lemma:basic-nxn} (\ref{basic-mu}), $|\mu| = 2m$ for some $m \in \{2, \ldots, 5\}$. Let $x \in \V(\Gamma)$. We apply (\ref{eq:b=c}) with $i = 2$ and $c_2(x,y) = |\mu| = 2m$ for all $y \in \Gamma_2(x)$ to count the number of edges between $\Gamma(x)$ and $\Gamma_2(x)$, and obtain
	\[ k_2(x) = \frac{1}{2m} \sum_{y \in \Gamma(x)} b_1(x,y) = \frac{1}{2m} \cdot 5^2(5-1)^2 = \frac{200}{m}. \]
So $m$ divides $200$, and thus $m \neq 3$. If $m = 5$ then Theorem \ref{theorem:srg} states that $\Gamma$ is strongly regular with parameters $(N, k, \lambda, \nu) = (66, 25, 8, 10)$. However
	\[ \frac{1}{2} \left( N - 1 \pm \frac{(N-1)(\nu-\lambda) - 2k}{\sqrt{(\nu-\lambda)^2 + 4(k-\nu)}} \right) = \frac{1}{2}(65 \pm 10), \]
neither of which is an integer, so there is no strongly regular graph having these parameters by \cite[Theorem 3.1]{srg}. Thus $m \neq 5$. Hence $m = 2$ or $4$, and $|\mu| = 4$ or $8$. If $|\mu| = 4$ then $\Gamma \cong J(10,5)$ or $\frac{1}{2}J(10,5)$ by \cite[Theorem 1]{4by4}. This proves \eqref{mu=8}.
\end{proof}

\begin{lemma} \label{lemma:5x5}
Assume that $\Gamma$ is connected and locally $5 \times 5$ grid, and that all $\mu$-graphs in $\Gamma$ have order at least $8$. For any $x \in \V(\Gamma)$:
	\begin{enumerate}[(1)]
	\item \label{5x5-clique} $\Gamma_2(x)$ does not contain any $6$-clique of $\Gamma$; and
	\item \label{eccentricity} $\epsilon(x) = 3$, where the eccentricity $\epsilon$ is as defined in (\ref{eq:epsilon}).
	\end{enumerate}
\end{lemma}

\begin{proof}
Suppose that $C \subseteq \Gamma_2(x)$ for some vertex $x$ and $6$-clique $C$. Then by Lemma \ref{lemma:cliques-bounded} (\ref{distance2-bounded}), all $y \in C$ satisfy $c_2(x,y) = 8$. Furthermore the six graphs $\mu(x,y)$, for $y \in C$, satisfy the conditions described in Lemma \ref{lemma:mu-clique}, namely, each pair of these six $\mu$-graphs of order $8$ is either disjoint or intersects in an edge (Lemma \ref{lemma:mu-clique} (\ref{mu-meet})) and the set of such edges forms a matching of $\left(\bigcup_{y \in C \cap \Gamma_2(x)} \mu(x,y)\right) \setminus C$ (Lemma \ref{lemma:mu-clique} (\ref{matching}); note that $C \cap \Gamma(x) = \varnothing$ since $C \subseteq \Gamma_2(x)$). However, a computer search using \textsc{Magma} \cite{magma} establishes that there is no set of subgraphs of $K_5 \square K_5$ that satisfy these conditions. (See Section \ref{sec:appendix} for the \textsc{Magma} code used.) Therefore $C \nsubseteq \Gamma_2(x)$. This proves statement (\ref{5x5-clique}).

To prove statement (\ref{eccentricity}), assume first that $\epsilon(x) = 2$ for some $x \in \V(\Gamma)$. Then $k_3(x) = 0$, so that for any $y \in \Gamma_2(x)$ and $6$-clique $C$ containing $y$, either $d_\Gamma(x,C) = 1$ or $C \subseteq \Gamma_2(x)$. But $C \nsubseteq \Gamma_2(x)$ by statement (\ref{5x5-clique}). So $d_\Gamma(x,C) = 1$ for any such $C$, and it follows that $c_2(x,y) = 10$ (for otherwise $c_2(x,y) = 8$, and so $d_\Gamma(x,C) = 2$ for some $6$-clique $C$ containing $y$ by Lemma \ref{lemma:maxcliques} (\ref{2m}), contradiction). Since $y$ is arbitrary, we then obtain $k_2(x) = 25(16)/10 = 40$ and
	\[ |\V(\Gamma)| = 1 + 25 + k_2(x) = 1 + 25 + 40 = 66. \]
If all vertices in $\Gamma$ have eccentricity $2$, then the above implies that $c_2(x,y) = 10$ for all $(x,y) \in \Gamma_2$. However this is impossible by Lemma \ref{lemma:mu-constant}. Thus $\epsilon(x') \geq 3$ for some $x' \in \V(\Gamma)$; since $\diam(\Gamma) \leq 3$ by Lemma \ref{lemma:mu-constant} \eqref{mu>8-lem}, we must then have $\epsilon(x') = 3$. In this case $k_3(x') \neq 0$, so by inequality (\ref{eq:k2,k3}) we have $k_3(x') = 1$ or $2$. Thus for some $y' \in \Gamma_2(x')$ and $6$-clique $C$ containing $y'$, we have $C \cap \Gamma_3(x') \neq \varnothing$, so that $d_\Gamma(x',C) = 2$. It follows from Lemma \ref{lemma:maxcliques} that $c_2(x',y') \neq 10$, so $c_2(x',y') = 8$. Hence $k_{2,8}(x') \neq 0$, where $k_{2,8}$ is as defined in (\ref{eq:k_2,2m}), and $\ell_{x'} := k_{2,8}(x')/5 \neq 0$. By equation (\ref{eq:k2-sum2}), $k_2(x') = 40 + \ell_{x'}$. Hence $k_2(x') > 40$, so that
	\[ |\V(\Gamma)| = 1 + 25 + k_2(x') + k_3(x') > 1 + 25 + 40 = 66, \]
contradiction. Therefore no vertex in $\Gamma$ has eccentricity $2$.
\end{proof}

For $x \in \V(\Gamma)$ and $m \in \{2, \ldots, n\}$ let
	\begin{equation} \label{eq:Gam_2,2m}
	\Gamma_{2,2m}(x) := \big\{ y \in \Gamma_2(x) \ : \ c_2(x,y) = 2m \big\}.
	\end{equation}

\begin{lemma} \label{lemma:5x5-order}
Assume that $\Gamma$ is locally $5 \times 5$ grid, and that all $\mu$-graphs in $\Gamma$ have order at least $8$. Then either:
	\begin{enumerate}[(1)]
	\item \label{order78} $|\V(\Gamma)| = 78$ and $c_2(x,y) = 8$ for all $x, y \in \Gamma$ with $d_\Gamma(x,y) = 2$, or
	\item \label{order72} $|\V(\Gamma)| = 72$, and with respect to any vertex $x$, $\Gamma$ has distance diagram as in Figure \ref{figure:diagram-5x5}.
	\end{enumerate}
\end{lemma}

\begin{figure}
\begin{center}
\begin{pspicture}(-5,-3)(5,1.5)
\rput(-5,0){\cnode*(0,0){3pt}{k0} \rput(-0.6,0){$x$}}
\rput(-1.9,0){\rput(0,0){\ovalnode{k1}{\parbox[c]{1.1cm}{\centering\small $\Gamma(x)$ \\ ($25$)}}} \rput(0.4,-0.85){\small $r$} \rput(1.1,0.25){\small $8$} \rput(0,1.3){\small $\Gamma_{2,8}(x')$} \rput{90}(0,0.9){\small $=$}}
\rput(0,-2.25){\rput(0,0){\ovalnode{k210}{\parbox[c]{1.2cm}{\centering\small $\Gamma_{2,10}(x)$ \\ $(20)$}}} \rput(-1,0.65){\small $10$} \rput(1,0.65){\small $10$} \rput(1.9,0){\small $= \Gamma_{2,10}(x')$}}
\rput(1.9,0){\rput(0,0){\ovalnode{k28}{\parbox[c]{1.1cm}{\centering\small $\Gamma_{2,8}(x)$ \\ ($25$)}}} \rput(0,1.3){\small $\Gamma(x')$} \rput{90}(0,0.9){\small $=$} \rput(-1.1,0.25){\small $8$} \rput(-0.4,-0.85){\small $s$}}
\rput(5,0){\cnode*(0,0){3pt}{k3} \rput(0.6,0){$x'$}}
\ncline{k0}{k1} \ncline{k1}{k28} \ncline{k28}{k3} \ncline{k1}{k210} \ncline{k28}{k210}
\end{pspicture}
\end{center}
\caption{Distance diagram for $\Gamma$ in the proof of Lemma \ref{lemma:5x5-order} (\ref{order72})}
\label{figure:diagram-5x5-proof}
\end{figure}
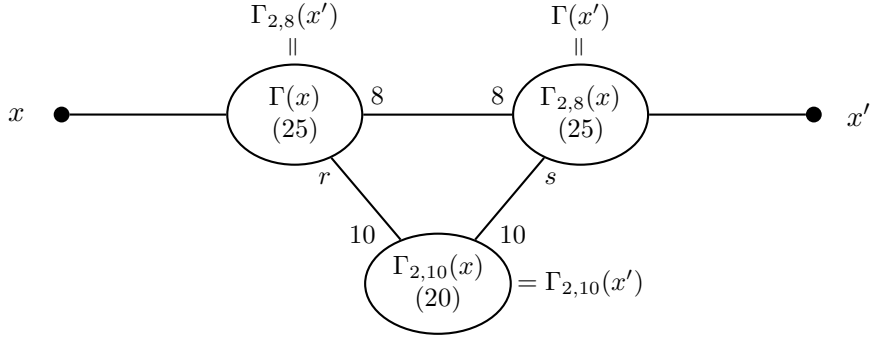

\begin{proof}
Let $x \in \V(\Gamma)$. Then $\epsilon(x) = 3$ by Lemma \ref{lemma:5x5} (\ref{eccentricity}), and hence $k_3(x) = 1$ or $2$ by the second inequality in (\ref{eq:k2,k3}). Also it follows from Lemma \ref{lemma:maxcliques} that $c_2(x,y) = 8$ for some $y \in \Gamma_2(x)$, and thus $k_{2,8}(x) \neq 0$. By (\ref{eq:k2-sum2}) we have $k_2(x) = 40 + \ell_x$ where $\ell_x := k_{2,8}(x)/5$, and $\ell_x \leq 5(4)/2 = 10$ by Lemma \ref{lemma:ell}. Hence
	\[ |\V(\Gamma)| = 1 + 25 + k_2(x) + k_3(x) = 1 + 25 + (40 + \ell_x) + k_3(x) = 66 + \ell_x + k_3(x). \]
Recall that $6$ divides $|\V(\Gamma)|$ by Lemma \ref{lemma:order-bound}. Hence $6$ divides $\ell_x + k_3(x)$, so the only possibilities for $(k_3(x),\ell_x)$ are $(1,5)$, $(2,4)$, and $(2,10)$. These yield $|\V(\Gamma)| = 72$ for $(k_3(x),\ell_x) \in \{ (1,5), (2,4) \}$, and $|\V(\Gamma)| = 78$ for $(k_3(x),\ell_x) = (2,10)$.

Assume that $|\V(\Gamma)| = 78$. It follows from the above that for any $x \in \V(\Gamma)$ we have $(k_3(x),\ell_x) = (2,10)$, so $k_2(x) = 40 + \ell_x = 50$ and $k_{2,8}(x) = 5\ell_x = 50$. Thus $k_{2,10}(x) = k_2(x) - k_{2,8}(x) = 0$, implying that $c_2(x,y) = 8$ for all $y \in \Gamma_2(x)$. Since $x$ is arbitrary, this means that $c_2$ is independent of $x$ or $y$, and thus all $\mu$-graphs in $\Gamma$ have order $8$. This proves statement (\ref{order78}).

For the remainder of the proof assume that $|\V(\Gamma)| = 72$. Then for any $x \in \V(\Gamma)$ we have $(k_3(x),\ell_x) \in \{(1,5), (2,4)\}$. In each case $\ell_x < 10$, and hence
	\[ k_{2,10}(x) = k_2(x) - k_{2,8}(x) = (40 + \ell_x) - 5\ell_x = 40 - 4\ell_x > 40 - 4(10) = 0. \]
So $\Gamma_{2,10}(x) \neq \varnothing$. By Lemma \ref{lemma:parameters} (\ref{b2}), $b_2(x,y) = 0$ for any $y \in \Gamma_{2,10}(x)$, and hence for any $x' \in \Gamma_3(x)$, $\Gamma(x') \subseteq \Gamma_{2,8}(x) \cup \Gamma_3(x)$. Thus $|\Gamma_{2,8}(x) \cup \Gamma_3(x)| \geq |\Gamma(x') \cup \{x'\}| = 26$. If $(k_3(x),\ell_x) = (2,4)$ for some $x \in \V(\Gamma)$, then $k_{2,8}(x) = 5\ell_x = 20$ and $|\Gamma_{2,8}(x) \cup \Gamma_3(x)| = k_{2,8}(x) + k_3(x) = 22$, contradiction. It follows that $(k_3(x),\ell_x) = (1,5)$ for all $x \in \V(\Gamma)$. Thus $k_{2,8}(x) = 5\ell_x = 25$, $k_{2,10}(x) = 40 - 4\ell_x = 20$, and there is a unique vertex $x' \in \Gamma_3(x)$. Hence, replacing $x$ by $x'$ in the above, $\Gamma(x') \subseteq \Gamma_{2,8}(x) \cup \Gamma_3(x) = \Gamma_{2,8}(x) \cup \{x'\}$. So $\Gamma(x') \subseteq \Gamma_{2,8}(x)$. Since $|\Gamma_{2,8}(x)| = k_{2,8}(x) = 25 = |\Gamma(x')|$, it follows that $\Gamma(x') = \Gamma_{2,8}(x)$, which in turn implies that $\Gamma_2(x') = \Gamma(x) \cup \Gamma_{2,10}(x)$.

Counting the number of edges between $\Gamma(x)$ and $\Gamma(x')$, and using the fact that $c_2(x,z) = 8$ for all $z \in \Gamma(x')$, we find that
	\[ \sum_{y \in \Gamma(x)} c_2(x',y) = \sum_{z \in \Gamma(x')} c_2(x,z) = 25(8). \]
Since $|\Gamma(x)| = 25$ and $c_2(x',y) \geq 8$ for all $y \in \Gamma(x)$, we must have $c_2(x',y) = 8$ for any $y \in \Gamma(x)$. Thus $\Gamma(x) \subseteq \Gamma_{2,8}(x')$. Note that also $(k_3(x'),\ell_{x'}) = (1,5)$, and thus $|\Gamma_{2,8}(x')| = k_{2,8}(x') = 5\ell_{x'} = 25 = |\Gamma(x)|$. Therefore $\Gamma_{2,8}(x') = \Gamma(x)$; since $\Gamma_2(x') = \Gamma_{2,8}(x') \cup \Gamma_{2,10}(x')$, this implies that $\Gamma_{2,10}(x') = \Gamma_{2,10}(x)$. This yields the distance diagram in Figure \ref{figure:diagram-5x5-proof}. Using the fact that $b_1 = 16$, we find that $r = s = b_1 - 8 = 8$.

Since $x$ is arbitrary, we then get the distance diagram in Figure \ref{figure:diagram-5x5-proof}, with the remaining parameters obtained using the fact that $\val(\Gamma) = 25$.
\end{proof}

\begin{proof}[Proof of Theorem \ref{maintheorem:5by5}]
Assume that $\Gamma$ is connected and locally $5 \times 5$ grid. It follows from Lemma \ref{lemma:basic-nxn} \eqref{basic-mu} that the $\mu$-graphs have order at most $10$. By Lemma \ref{lemma:mu-constant} \eqref{mu=8}, if all $\mu$-graphs have the same order then either $|\mu| = 4$ or $|\mu| = 8$, and if $|\mu| = 4$ then $\Gamma \cong J(10,5)$ or $\frac{1}{2} J(10,5)$. So Theorem \ref{maintheorem:5by5} \eqref{mu-constant} holds for the case $|\mu| = 4$. Moreover, it follows that not all $\mu$-graphs can have order $10$, so there exists $(x,y) \in \Gamma_2$ such that $\mu(x,y)$ has order $c_2(x,y) \leq 8$.

Thus we may assume that all $\mu$-graphs have order at least $8$. Then $\diam(\Gamma) \leq 3$ by Lemma \ref{lemma:mu-constant} \eqref{mu>8-lem}, and $|\V(\Gamma)| = 72$ or $78$ by Lemma \ref{lemma:5x5-order}. If $|\V(\Gamma)| = 78$ then by Lemma \ref{lemma:5x5-order} \eqref{order78} all $\mu$-graphs have order equal to $8$, and by Theorem \ref{maintheorem:order-diameter} \eqref{=2(n-1)} the graph $\Gamma$ is a distance-regular antipodal triple cover of $K_{26}$ with intersection array $(25, 16, 1; 1, 8, 25)$. This proves Theorem \ref{maintheorem:5by5} \eqref{mu>8} (ii). This also completes the proof of Theorem \ref{maintheorem:5by5} \eqref{mu-constant}.

To complete the proof of Theorem \ref{maintheorem:5by5} \eqref{mu>8} (i) we need to consider the case $|\V(\Gamma)| =  72$. Then by Lemma \ref{lemma:5x5-order} \eqref{order72}, with respect to any $x \in \V(\Gamma)$, $\Gamma$ has distance diagram as in Figure \ref{figure:diagram-5x5}. Let $y \in \Gamma(x)$, let $x'$ be the unique vertex in $\Gamma_3(x)$, and let $y'$ be the unique vertex in $\Gamma_3(y)$. Clearly $y' \in \Gamma(x)$, since otherwise $x$ is a common neighbour of $y$ and $y'$. Also $y' \neq x'$, since $d_\Gamma(y,x') = 2$. Suppose that $y' \in \Gamma_{2,10}(x)$. Since $\Gamma_{2,10}(x) = \Gamma_{2,10}(x')$, we have $c_2(x',y') = 10$. So $x' \in \Gamma_{2,10}(y')$. Applying Lemma \ref{lemma:5x5-order} (\ref{order72}) using $y$ in the place of $x$, we find that $\Gamma_{2,10}(y') = \Gamma_{2,10}(y)$, so $c_2(x',y) = 10$. But $y \in \Gamma(x) = \Gamma_{2,8}(x')$, so $c_2(x',y) = 8$, contradiction. Therefore $y' \notin \Gamma_{2,10}(x)$, and so $y' \in \Gamma_{2,8}(x) = \Gamma(x')$. Consequently, for any $z \in \Gamma_{2,10}(x)$, the unique $z' \in \Gamma_3(x)$ also lies in $\Gamma_{2,10}(x)$. Moreover, for any $x, x', y, y' \in \V(\Gamma)$ with $d_\Gamma(x,x') = d_\Gamma(y,y') = 3$, we have $x \sim_\Gamma y$ if and only if $x' \sim_\Gamma y'$.

Now consider the quotient graph $\Gamma_\mathcal{P}$ of $\Gamma$ with respect to the partition $\mathcal{P} = \big\{ \{x,x'\} \ : \ (x,x') \in \Gamma_3 \big\}$. The graph $\Gamma_\mathcal{P}$ has diameter $2$ and valency $25$. Let $x \in \V(\Gamma)$. Since each $y \in \Gamma(x)$ has $8$ neighbours in $\Gamma(x)$ and $8$ neighbours in $\Gamma_{2,8}(x)$, the vertex $\{y,y'\}$ of $\Gamma_\mathcal{P}$ has $16$ neighbours in $\Gamma_\mathcal{P}(\{x,x'\})$. Any $z \in \Gamma_{2,10}(x)$ has $10$ neighbours in $\Gamma(x)$ and $10$ neighbours in $\Gamma_{2,8}(x)$; clearly if $u \in \Gamma(z) \cap \Gamma(x)$ and $v \in \Gamma(z) \cap \Gamma_{2,8}(x)$ then $u$ and $v$ belong to different blocks of $\mathcal{P}$. So $\{z,z'\}$ has $20$ neighbours in $\Gamma_\mathcal{P}(\{x,x'\})$. Thus $\Gamma_\mathcal{P}$ is strongly regular with parameters $(36,25,16,20)$. There is a unique strongly regular graph with these parameters, namely, the complement $K_6 \times K_6$ of the $6 \times 6$ grid $K_6 \square K_6$ (\cite{srgdatabase} or \cite[Table 1.1, line 13]{srgBM}).

This completes the proof of Theorem \ref{maintheorem:5by5}.
\end{proof}

%% file: locallygrid-sec7-appendix.tex
\section{Appendix: Magma program for Lemma \ref{lemma:5x5} (\ref{5x5-clique})} \label{sec:appendix}

Assume that $\Gamma$ is locally $5 \times 5$ grid, and that all $\mu$-graphs of $\Gamma$ have order at least $8$. Let $x \in \V(\Gamma)$. We want to determine if some $6$-clique $C$ is contained in $\Gamma_2(x)$. If such a clique exists, then each $\mu$-graph $\mu(x,y)$, for $y \in C$, is an induced subgraph of $K_5 \square K_5$ and is either an $8$-cycle or a disjoint union of two $4$-cycles, and the set of these six $\mu$-graphs satisfies the conditions in Lemma \ref{lemma:mu-clique}, namely:
	\begin{enumerate}
	\item If $S$ is union of vertex sets of these six $\mu$-graphs, then each vertex in $S$ lies in exactly two $\mu$-graphs.
	\item Any two distinct $\mu$-graphs are either disjoint or have exactly one common edge.
	\item The set of edges which lie in two $\mu$-graphs (as in 2.) form a perfect matching of $S$ (since $C \cap \Gamma(x) = \varnothing$).
	\end{enumerate}
For the computation we replaced condition (1) with the following weaker condition:
	\begin{enumerate}
	\item[(1')] Each vertex in S lies in at most two $\mu$-graphs.
	\end{enumerate}

We denoted by $\texttt{Cyc8}$ and $\texttt{Cyc44}$, respectively, the set of all induced subgraphs of $K_5 \square K_5$ which are $8$-cycles, and the set of all induced subgraphs of $K_5 \square K_5$ which are unions of two disjoint $4$-cycles. Since $\Aut(K_5 \square K_5)$ is transitive on each of the sets $\texttt{Cyc8}$ and $\texttt{Cyc44}$, so we assumed without loss of generality that:
	\begin{enumerate}
	\item[(4)] One of the six graphs is a fixed graph $\texttt{mu}$.
	\end{enumerate}

We considered two cases, one with $\texttt{mu} \in \texttt{Cyc8}$ and the other with $\texttt{mu} \in \texttt{Cyc44}$. For each of these cases we used \textsc{Magma} to enumerate all sets consisting of $6$ induced subgraphs of $K_5 \square K_5$ satisfying the conditions (1'), (2), and (4). In each case we managed to find as many as five such subgraphs but there were no sets of six.

The following is our \text{Magma} code.

\subsection*{Step 1.} We constructed the graph $K_5 \square K_5$ and the sets $\texttt{Cyc8}$ and $\texttt{Cyc44}$.

\bigskip \noindent \tt
n := 5; \\[6pt]
vertices := \{ <a,b> : a,b in [1..n] \}; \\
edges := \{\{u,v\} : u,v in vertices | u ne v and (u[1] eq v[1] or u[2] eq v[2])\}; \\
grid,V,E := Graph< vertices | edges >; \\[6pt]
Cyc8 := \{ X : X in Subsets(Set(V),2*(n-1)) | IsIsomorphic(sub< grid | X >, \\ PolygonGraph(2*(n-1))) \}; \\[6pt]
Cyc44 := \{ X : X in Subsets(Set(V),2*(n-1)) | IsIsomorphic(sub< grid | X >, \\ Union(PolygonGraph(4),PolygonGraph(4))) \};

\subsection*{Step 2.} \rm We constructed the set of all $2$-sets of graphs in $\texttt{Cyc8} \cup \texttt{Cyc44}$ which satisfy (2):

\bigskip \noindent \tt
U := \{ \{@ X1,X2 @\} : X1,X2 in Cyc44 join Cyc8 | IsDisjoint(X1,X2) or ( \#(X1 \\ meet X2) eq 2 and IsIsomorphic(sub< grid | X1 meet X2 >, CompleteGraph(2)) ) \};

\subsection*{Step 3.} \rm We constructed the fixed graph $\texttt{mu} \in \texttt{Cyc8}$:

\bigskip \noindent \tt
mu := \{V!<1,1>, V!<1,2>, V!<2,2>, V!<2,3>, V!<3,3>, V!<3,4>, V!<4,4>, V!<4,1>\};

\subsection*{Step 4.} \rm We constructed the sets $\texttt{W}$, $\texttt{X}$, $\texttt{Y}$, and $\texttt{Z}$ of all $3$-, $4$-, $5$-, and $6$-sets, respectively, of $\mu$-graphs satisfying (1'), (2), and (4). Note that \textsc{Magma} returned non-empty sets $\texttt{W}$, $\texttt{X}$, and $\texttt{Y}$, but that $\texttt{Z}$ is empty.

\bigskip \noindent \tt
W := \{ \{@ mu, x[1], x[2] @\} : x in U | mu notin x and forall(i)\{ i : i in \\ {[1..\#x]} | \{@ x[i], mu @\} in U \} and IsDisjoint(mu, x[1] meet x[2]) \}; \\[6pt]
X := \{ \{@ x[1], x[2], x[3], z @\} : x in W, z in Cyc44 join Cyc8 | z notin x \\ and forall(i)\{ i : i in [1..\#x] | \{@ x[i], z @\} in U \} and forall(i)\{ \{i,j\} : \\ i,j in [1..\#x] | i eq j or IsDisjoint(z, x[i] meet x[j]) \} \}; \\[6pt]
Y := \{ \{@ x[1], x[2], x[3], x[4], z @\} : x in X, z in Cyc44 join Cyc8 | z \\ notin x and forall(i)\{ i : i in [1..\#x] | \{@ x[i], z @\} in U \} and forall(i)\{ \\ \{i,j\} : i,j in [1..\#x] | i eq j or IsDisjoint(z, x[i] meet x[j]) \} \}; \\[6pt]
Z := \{ \{@ x[1], x[2], x[3], x[4], x[5], z @\} : x in Y, z in Cyc44 join Cyc8 | \\ z notin x and forall(i)\{ i : i in [1..\#x] | \{@ x[i], z @\} in U \} and \\ forall(i)\{ \{i,j\} : i,j in [1..\#x] | i eq j or IsDisjoint(z, x[i] meet x[j]) \} \};


\subsection*{Step 5} \rm Finally we repeated steps 3 to 4 for the fixed graph $\texttt{mu} \in \texttt{Cyc44}$, namely

\bigskip \noindent \tt
mu := \{V!<1,1>, V!<1,2>, V!<2,2>, V!<2,1>, V!<3,3>, V!<3,4>, V!<4,4>, V!<4,3>\};

\bigskip \noindent \rm
and again the collection \texttt{Z} of $6$-sets was empty.


\rm